\documentclass[11pt]{amsart}
\usepackage{amsmath, amssymb}
\usepackage[all]{xy}
\usepackage{verbatim}
\usepackage{multirow}
\usepackage{enumerate}

\newcommand{\qqed}{\begin{flushright}
$\Box$
\end{flushright}}
\newcommand{\R}{\mathbb{R}}
\newcommand{\C}{\mathbb{C}}

\newtheorem{thm}{Theorem}[section]
\newtheorem{lem}[thm]{Lemma}
\newtheorem{prop}[thm]{Proposition}
\newtheorem{defn}[thm]{Definition}
\newtheorem{cor}[thm]{Corollary}

\DeclareMathOperator{\range}{range}
\DeclareMathOperator{\spann}{span}
\DeclareMathOperator{\rank}{rank}

\newcommand{\mA}{\mathfrak{A}(L^2(0, \infty))_* ^+}

\DeclareMathOperator{\nullspace}{nullspace}
\author{Christopher Jankowski} \thanks{Supported by a Graduate Student Fellowship at the
University of Pennsylvania and later by the Skirball Foundation via
the Center for Advanced Studies in Mathematics at Ben-Gurion
University of the Negev.}
\address{ \newline Department of Mathematics \newline
Ben-Gurion University of the Negev\newline P.O. Box 653 \newline
Be'er Sheva 84105, Israel} \email{cjankows@math.bgu.ac.il}
\title{Unital $q$-positive maps on $M_2(\C)$ and a related
$E_0$-semigroup result}
\begin{document}
\begin{abstract}
From previous work, we know how to obtain type II$_0$ $E_0$-semigroups
using boundary weight doubles $(\phi, \nu)$, where $\phi: M_n(\C)
\rightarrow M_n(\C)$ is a unital $q$-positive map and $\nu$ is a
normalized unbounded boundary weight over $L^2(0, \infty)$.
In this paper, we classify the unital $q$-positive maps $\phi:
M_2(\C) \rightarrow M_2(\C)$.  We find that every unital $q$-pure
map $\phi: M_2(\C) \rightarrow M_2(\C)$ is either rank one or
invertible. We also examine the case $n=3$, finding the limit maps
$L_\phi$ for all unital $q$-positive maps $\phi: M_3(\C) \rightarrow
M_3(\C)$. In conclusion, we present a cocycle conjugacy result for
$E_0$-semigroups induced by boundary weight doubles $(\phi, \nu)$
when $\nu$ has the form $\nu(\sqrt{I - \Lambda(1)} B \sqrt{I -
\Lambda(1)})=(f,Bf).$
\end{abstract}
\maketitle
\section{Introduction}
A linear map $\phi: M_n(\C) \rightarrow M_n(\C)$ with no negative
eigenvalues is said to be $q$-positive if $\phi(I + t \phi)^{-1}$ is
completely positive for all $t \geq 0$.  This class of maps has
recently played a key role in constructing $E_0$-semigroups in
\cite{Me}.  Let $H$ be a separable Hilbert space whose inner product
$( \ , \ )$ is conjugate-linear in its first entry and linear in its
second.  An $E_0$-semigroup $\alpha=\{\alpha_t\}_{t\geq 0}$ is a
weakly continuous semigroup of unital $*$-endomorphisms of $B(H)$.
Every $E_0$-semigroup $\alpha$ is assigned one of three types based
on intertwining semigroups called units. A \textit{unit} for
$\alpha$ is a strongly continuous semigroup $V=\{V_t\}_{t \geq 0}$
of operators in $B(H)$ such that $\alpha_t(A)V_t=V_tA$ for all $t
\geq 0$ and $A \in B(H)$. Let $\mathcal{U}_\alpha$ be the set of
units for $\alpha$.  If $\mathcal{U}_\alpha$ is nonempty, we say
$\alpha$ is spatial. If, for all $t \geq 0$, the closed linear span
of the set $\{ U_{1}(t_1) \cdots U_n(t_n) f: f \in H, t_i \geq 0
\textrm{ and } U_i \in \mathcal{U}_\alpha \ \forall \ i, \sum t_i =
t \}$ is $H$, we say $\alpha$ is completely spatial.  If $\alpha$ is
completely spatial, we say $\alpha$ is of type I, while if $\alpha$
is spatial but is not completely spatial, we say $\alpha$ is of type
II.  If $\alpha$ has no units, we say $\alpha$ is of type III.  Each
spatial $E_0$-semigroup is given an index $n \in \mathbb{Z}_{\geq 0}
\cup \{\infty\}$ which depends on the structure of its units and is
invariant under cocycle conjugacy.

We can naturally construct $E_0$-semigroups over symmetric and
antisymmetric Fock spaces using the right shift semigroup on $K
\otimes L^2(0, \infty)$, obtaining the CCR and CAR flows of rank
$dim(K)$. These yield all non-trivial type I $E_0$-semigroups in
terms of cocycle conjugacy: If $\alpha$ is of type I$_n$ (type I,
index $n$) for $n \in \mathbb{N} \cup \{\infty\}$, then $\alpha$ is
cocycle conjugate to the CCR flow of rank $n$ (see \cite{arvindex}).
The classification of $E_0$-semigroups of types II and III is far
more complicated, however. Uncountably many examples of both types
are known and have been exhibited through greatly differing methods
(see, for example, \cite{izumi}, \cite{bigpaper}, and  \cite{T2}).
Using Bhat's dilation theorem (\cite {Bhat}), Powers showed in
\cite{hugepaper} that every spatial $E_0$-semigroup is induced by
the boundary weight map of a $CP$-flow over $K \otimes L^2(0,
\infty)$ for a separable Hilbert space $K$.  He investigated the
case when $dim(K)=1$ in \cite{bigpaper}, exhibiting uncountably many
mutually non-cocycle conjugate type II$_0$ $E_0$-semigroups using
boundary weights over $L^2(0, \infty)$.  He also began to explore
the case when $K$ is $2$-dimensional by combining Schur maps with
boundary weights.  This approach was generalized to the case when
$1<dim(K)<\infty$ in \cite{Me}, where the theory of boundary weight
doubles was introduced.

A boundary weight double is a pair $(\phi, \nu)$, where $\phi:
M_n(\C) \rightarrow M_n(\C)$ is a $q$-positive map and $\nu$ is a
positive boundary weight over $L^2(0, \infty)$ (we write $\nu \in
\mA$). If $\phi$ is unital and $\nu$ is normalized and unbounded (in
which case we call $\nu$ a type II Powers weight), then $(\phi,
\nu)$ induces a unital $CP$-flow over $\C^n$ whose Bhat minimal
dilation is a type II$_0$ $E_0$-semigroup. Comparing
$E_0$-semigroups induced by boundary weight doubles in terms of
cocycle conjugacy becomes easier if we focus on the \textit{$q$-pure
maps}, which are $q$-positive maps with the smallest possible
structure of $q$-subordinates (see Definition \ref{qpure}).  The
unital $q$-pure maps which are either rank one or invertible have
all been classified in \cite{Me}:  The unital rank one $q$-pure maps
$\phi: M_n(\C) \rightarrow M_n(\C)$ are implemented by faithful
states in $M_n(\C)^*$, while the unital invertible $q$-pure maps are
a particular class of Schur maps (see Theorems \ref{statesbig} and
\ref{inverts} for a summary).

Our main goal in this paper is to begin the general classification
of all unital $q$-positive maps $\phi: M_n(\C) \rightarrow M_n(\C)$,
with the particular aim of finding all such maps which are $q$-pure.
Our second goal is to prove cocycle conjugacy comparison results for
boundary weight doubles $(\phi, \nu)$ and $(\psi, \nu)$
when $\phi$ and $\psi$ are not $q$-pure.  We should note that we are only interested in
identifying a $q$-positive map $\phi$ up to a particular notion of
equivalence which we call conjugacy.  More precisely, for each
$q$-positive $\phi: M_n(\C) \rightarrow M_n(\C)$ and unitary $U \in
M_n(\C)$, we can form a new $q$-positive map $\phi_U: M_n(\C)
\rightarrow M_n(\C)$ by defining $\phi_U(A)=U^*\phi(UAU^*)U$ for all
$A \in M_n(\C)$.  If $\phi$ is unital and $\nu$ is a type II Powers
weight of the form $\nu(\sqrt{I -
\Lambda(1)} B \sqrt{I-\Lambda(1)})=(f,Bf)$, then $(\phi, \nu)$ and
$(\phi_U, \nu)$ induce cocycle conjugate $E_0$-semigroups
(Proposition \ref{dumb}).  In fact, a much stronger result holds:  If
$\phi$ is unital and $\nu$ is any type II Powers weight, then
$(\phi, \nu)$ and $(\phi_U, \nu)$ induce conjugate $E_0$-semigroups
(\cite{DJ}). Motivated by this fact, we say that $q$-positive maps
$\phi, \psi: M_n(\C) \rightarrow M_n(\C)$ are \textit{conjugate} if
$\psi=\phi_U$ for some unitary $U \in M_n(\C)$.

Let $\mathcal{E}_n$ be the set of all unital completely positive
maps $\Phi: M_n(\C) \rightarrow M_n(\C)$ such that $\Phi^2=\Phi$.
This is merely the set of all limits $L_\phi = \lim_{t \rightarrow
\infty} t \phi(I + t \phi)^{-1}$ for unital $q$-positive maps $\phi:
M_n(\C) \rightarrow M_n(\C)$.  This limiting method has already
appeared in \cite{Me}, where it was vital in classifying the unital
rank one $q$-pure maps on $M_n(\C)$. We find all elements of
$\mathcal{E}_2$ and $\mathcal{E}_3$ up to conjugacy.  Using this
result, we classify the unital $q$-positive maps $\phi: M_2(\C)
\rightarrow M_2(\C)$, finding that there is no unital $q$-positive
map $\phi: M_2(\C) \rightarrow M_2(\C)$ of rank $3$ (Proposition
\ref{E2}). Moreover, we find that that the only unital $q$-pure maps
$\phi: M_2(\C) \rightarrow M_2(\C)$ are either rank one or
invertible (Theorem \ref{qpurem2}).  We also show that any unital
$q$-positive map $\phi: M_3(\C) \rightarrow M_3(\C)$ which
annihilates a nonzero positive matrix cannot be $q$-pure (see
Proposition \ref{notqpure}).

In conclusion, we compare $E_0$-semigroups
formed by boundary weight doubles 
$(\phi, \nu)$ and $(\psi, \nu)$ in the case that
$\phi: M_n(\C) \rightarrow M_n(\C)$ ($n \geq 2$) is any unital
rank one $q$-positive map, $\psi: M_k(\C) \rightarrow M_k(\C)$
is any unital $q$-positive map such that $L_\psi$ is a Schur map,
and $\nu$ is a type II Powers weight of the form
$\nu(\sqrt{I -
\Lambda(1)} B \sqrt{I-\Lambda(1)})=(f,Bf)$
(Theorem \ref{genl}).  This
result substantially generalizes a consequence of Theorems 5.4 and 6.12
of \cite{Me}.  

\section{Background}
\subsection{Completely positive and $q$-positive maps}

Let $\phi: B(K) \rightarrow B(H)$ be a linear map.  We say that
$\phi$ is \textit{unital} if $\phi(I_K)=I_H$ and \textit{positive}
if $\phi(A)$ is positive whenever $A \in B(K)$ is positive. For each
$n \in \mathbb{N}$, define $\phi_n: M_n(B(K)) \rightarrow M_n(B(H))$
by

\begin{displaymath} \phi_n \left(\begin{array}{ccc} A_{11} & \cdots & A_{1n}
\\ \vdots & \ddots & \vdots \\ A_{n1} & \cdots & A_{nn}  \end{array}\right)
= \left(\begin{array}{ccc} \phi(A_{11}) & \cdots & \phi(A_{1n})
\\ \vdots & \ddots & \vdots \\ \phi(A_{n1}) & \cdots & \phi(A_{nn})
\end{array} \right).
\end{displaymath}
We say that $\phi$ is completely positive if $\phi_n$ is positive
for all $n \in \mathbb{N}$.  If $\phi$ is completely positive, then
$||\phi||=||\phi(I_K)||$.

We know from a result of Choi (see \cite{choi}) that a linear map
$\phi: M_n(\C) \rightarrow M_n(\C)$ is completely positive if and only if it
can be written in the form
$$\phi(A) = \sum_{i=1}^k S_i A S_i^*$$ for some integer $k \leq n^2$
and linearly independent $n \times n $ matrices $\{S_i\}_{i=1}^k$.
This result generalizes to normal completely positive maps between
$B(K)$ and $B(H)$ for separable Hilbert spaces $K$ and $H$ (see
\cite{arveson}).  Denote by $\{e_{ij}\}_{i,j=1}^n$ the set of
standard matrix units for $M_n(\C)$. Given any $M=\sum_{i,j=1}^n
a_{ij}e_{ij} \in M_n(\C)$, we can form a linear map $\phi: M_n(\C)
\rightarrow M_n(\C)$ by defining $\phi(A)= \sum_{i,j}
m_{ij}a_{ij}e_{ij}$ for all $A=\sum_{i,j=1}^n a_{ij}e_{ij} \in
M_n(\C)$. We call this the Schur map corresponding to $M$, and
denote it by the notation $\phi(A)= M \bullet A$.  We will
frequently use the fact that $\phi$ is completely positive if and
only if $M$ is positive (for a proof, see \cite{paulsen}).  By a
positive matrix we mean a self-adjoint matrix whose eigenvalues are
all nonnegative.


The construction of $E_0$-semigroups in \cite{Me} (as we will see in
Proposition \ref{bdryweight}) required a particular kind of completely
positive map:
\begin{defn}\label{qpos}  A linear map $\phi: M_n(\C) \rightarrow
M_n(\C)$ is
\emph{$q$-positive} if $\phi$ has no negative eigenvalues and
$\phi(I + t \phi)^{-1}$ is completely
positive for all $t \geq 0$.
\end{defn}

The condition that a completely positive map $\phi$ must have no
negative eigenvalues
in order to be $q$-positive
is certainly non-trivial,
as completely positive maps with negative eigenvalues exist in abundance.  One such
example is the Schur map $\phi: M_2(\C) \rightarrow M_2(\C)$ defined by
\begin{displaymath}
\phi \left(\begin{array}{cc} a_{11} & a_{12} \\ a_{21} & a_{22} \end{array}\right)
= \left(\begin{array}{cc} a_{11} & - a_{12} \\ - a_{21} & a_{22} \end{array} \right).
\end{displaymath}

Furthermore, even if $\phi$ is a completely positive map with no negative 
eigenvalues,
it does not necessarily follow that $\phi(I + t \phi)^{-1}$ is completely positive for
all $t \geq 0$.  In fact, for each $s \geq 0$, we can construct a completely
positive map $\phi$ which is not $q$-positive but which still satisfies the condition
that $\phi(I + t \phi)^{-1}$ is completely positive for all $0 \leq t \leq s$.  For this, let
$r \in (1, \sqrt{2}]$ and define a Schur map $\phi_r: M_2(\C) \rightarrow M_2(\C)$ by
\begin{displaymath}
\phi_r \left(\begin{array}{cc} a_{11} &  a_{12} \\ a_{21} & a_{22} \end{array} \right) =
\left(\begin{array}{cc} a_{11} & \frac{r(1+i) a_{12}}{2} \\ \frac{r(1-i)a_{21}}{2} & a_{22} \end{array}\right).
\end{displaymath}

In other words, $\phi_r(A) = M \bullet A$ for the positive matrix
\begin{displaymath}
M= \left(\begin{array}{cc} 1 & \frac{r(1+i)}{2} \\ \frac{r(1-i)}{2} & 1\end{array}\right).
\end{displaymath}
We find that
$\phi_r(I + t \phi_r)^{-1}(A) = M_{t} \bullet A$ for all $A \in M_n(\C)$ and $t \geq 0$, where
\begin{displaymath}
M_t = \left(\begin{array}{cc} \frac{1}{1+t} & \frac{r(1+i)}{2+ tr(1+i)} \\ \frac{r(1-i)}{2+tr(1-i)} &
\frac{1}{1+t}\end{array}\right).
\end{displaymath}
As noted previously, $A \rightarrow M_t \bullet A$ is completely positive
if and only if $M_t$ is a positive matrix.
Let $\lambda_1$ and $\lambda_2$ be the eigenvalues of $M_t$.  Since $\lambda_1+\lambda_2 =tr(M_t)>0$
and $\lambda_1 \lambda_2 = \det(M_t)$, $M_t$ is positive if and only if its determinant is
nonnegative.  A calculation
shows that for any given
$t \geq 0$, $\det(M_t)$ is
nonnegative if and only if $t \leq \frac{2-r^2}{2r(r-1)}.$
Therefore,
$\phi_r(I + t \phi_r)^{-1}$ ($t\geq 0$) is completely positive if and only if
$$t \leq \frac{2-r^2}{2r(r-1)}.$$
Let $s\geq 0$.  The values $(2-r^2)/(2r^2-2r)$ for $r \in (1, \sqrt{2}]$
clearly span $[0, \infty)$, so
$s= (2-r_0^2)/(2r_0^2-2r_0)$ for some $r_0 \in (1, \sqrt{2}]$. By the previous paragraph,
$\phi_{r_0}(I + t \phi_{r_0})^{-1}$ is
completely positive if $0 \leq t \leq s$ but is not completely positive if $t>s$.  This example
demonstrates that we cannot generally conclude that a map $\phi$ is $q$-positive if
$\phi(I + t \phi)^{-1}$ is completely positive
for all $t$ in some finite interval $J \subset \R_{\geq 0}$, no matter how large $J$ is.

There is a natural order structure for $q$-positive maps.
If $\phi, \psi:
M_n(\C) \rightarrow M_n(\C)$ are $q$-positive, we say that $\phi$ $q$-dominates
$\psi$ (i.e. $\phi \geq_q \psi$) if $\phi(I + t \phi)^{-1} - \psi(I + t \psi)^{-1}$
is completely positive for all $t \geq 0$.  As it turns out, for 
every $s \geq 0$, the map $\phi(I + s \phi)^{-1}$ is $q$-positive and
$\phi \geq_q \phi(I + s \phi)^{-1}$ (Proposition 4.1 of \cite{Me}).

\begin{defn}\label{qpure}  A $q$-positive map $\phi: M_n(\C) \rightarrow M_n(\C)$ is \emph{$q$-pure}
if its set of $q$-subordinates is $\{\phi(I + s \phi)^{-1}\}_{s \geq 0} \cup \{0\}$.
\end{defn}

\subsection{$E_0$-semigroups and $CP$-flows}
A result of Wigner in \cite{wigner} shows that every one-parameter
group $\alpha = \{\alpha_t\}_{t \in \R}$
of $*$-automorphisms of $B(H)$ is implemented by a strongly
continuous unitary group
$U=\{U_t\}_{t \in \R}$
in the sense that
$$\alpha_t(A)= U_t A U_t^*$$ for all $A \in B(H)$ and $t \geq 0$.
This leads us to ask how to characterize all suitable semigroups
of $*$-endomorphisms
of $B(H)$:
\begin{defn}
We say a family $\{\alpha_t\}_{t \geq 0}$ of $*$-endomorphisms of
$B(H)$ is an \emph{$E_0$-semigroup} if:
\begin{enumerate}[(i)]
\item $\alpha_{s+t}=\alpha_s \circ \alpha_t$ for all $s, t \geq 0$,
and $\alpha_0 (A)=A$ for all $A \in B(H)$.
\item  For each $f, g \in H$ and $A \in B(H)$, the inner product
$( f, \alpha_t (A) g)$ is continuous in $t$.
\item $\alpha_t (I)=I$ for all $t \geq 0$ (in other words, $\alpha$ is
\emph{unital}).
\end{enumerate}
\end{defn}

There are two different conditions under which we think of
$E_0$-semigroups as equivalent.  The first, and stronger condition,
is conjugacy, while the second condition, cocycle conjugacy, will be
our main focus in comparing $E_0$-semigroups.

\begin{defn}  Let $\alpha$ and $\beta$ be $E_0$-semigroups on $B(H_1)$ and
$B(H_2)$, respectively.  We say that $\alpha$ and $\beta$ are
\emph{conjugate} if there is a $*$-isomorphism $\theta$ from
$B(H_1)$ onto $B(H_2)$ such that $\theta \circ \alpha_t \circ \theta^{-1}= \beta_t
$ for all $t \geq 0$.

We say that $\alpha$ and $\beta$
are \emph{cocycle conjugate} if $\alpha$ is conjugate to $\beta'$,
where $\beta'$ is an $E_0$-semigroup on $B(H_2)$ satisfying the
following condition: For some strongly continuous family of
unitaries $U=\{U_t: t \geq 0\}$ acting on $H_2$ and satisfying
$U_{t+s}=U_t \beta_t(U_s)$ for all $s, t \geq 0$, we have $\beta'_t
(A) = U_t \beta_t (A) U_t^*$ for all $A \in B(H_2)$ and $t \geq 0$.
\end{defn}

Bhat's dilation theorem from \cite{Bhat} shows that we can obtain
$E_0$-semigroups from much more general semigroups of completely
positive maps called $CP$-semigroups.  A $CP$-flow is a
$CP$-semigroup acting on $B(K \otimes L^2(0, \infty))$ which is
intertwined by the right shift semigroup.  More specifically:

\begin{defn}  Let $H= K \otimes L^2(0, \infty)$, which we identify
with the space of $K$-valued measurable functions defined on $(0,
\infty)$ which are squre integrable.  Denote by $U=\{U_t\}_{t \geq
0}$ the right shift semigroup on $H$, so for all $f \in H$, $x \in
(0, \infty)$, and $t \geq 0$, we have $(U_tf)(x) = f(x-t)$ if $x>t$
and $(U_tf)(x)=0$ otherwise.

A strongly continuous semigroup $\alpha=\{\alpha_t: t \geq 0\}$ of
completely positive contractions of $B(H)$ into itself is called a
\textit{CP-flow} if $\alpha_t(A)U_t = U_t A$ for all $t \geq 0$ and
$A \in B(H)$.
\end{defn}

Unless otherwise specified, we will henceforth write $\{U_t\}_{t \geq 0}$
for the right shift semigroup acting on $K \otimes L^2(0, \infty)$.
Special functionals called boundary weights play an important role
in constructing $CP$-flows (see Definition 1.10 of \cite{markie} for
a more general definition and a detailed discussion):
\begin{defn}
Let $H= K \otimes L^2(0,\infty)$ and define $\Lambda: B(K)
\rightarrow B(H)$ by $$(\Lambda(A)f)(x)=e^{-x}Af(x)$$ for all
$A \in B(K)$, $f \in H$, and $x \in (0, \infty)$.  We denote
by $\mathfrak{A}(H)$ the linear space
$$\mathfrak{A}(H)= \sqrt{I - \Lambda(I_K)}B(H)\sqrt{I - \Lambda(I_K)}$$
and by $\mathfrak{A}(H)_*$ the linear functionals $\rho$ on $\mathfrak{A}$
of the form
$$\rho \Big(\sqrt{I - \Lambda(I_K)} A \sqrt{I - \Lambda(I_K)}\Big) = \eta(A)$$
for $A \in B(H)$ and $\eta \in B(H)_*$.  We call such functionals boundary
weights. 
\end{defn}

We can associate to every $CP$-flow $\alpha$ a boundary weight map
$\rho \rightarrow \omega(\rho)$ from $B(K)_*$ to $\mathfrak{A}(H)_*$ which
is related to $\alpha$ in the following manner.
Let $R_\alpha$
be the resolvent
$$R_\alpha(A) = \int_0 ^\infty e^{-t} \alpha_t(A) dt$$ of $\alpha$, and define
$\Gamma: B(H) \rightarrow B(H)$ by $\Gamma(A) = \int_0 ^\infty
e^{-t} U_t A U_t^* dt$ for all $A \in B(H)$.  Using hats to denote the predual
mappings, we have
$$\hat{R}_\alpha(\tau) = \hat{\Gamma}\Big(\omega(\hat{\Lambda}\tau) + \tau\Big)$$ for all
$\tau \in B(H)_*$.
If we let $\rho \rightarrow \omega_t(\rho)$ be the truncated boundary weight maps
\begin{eqnarray}\label{truncate} \omega_t(\rho)(A)=
\omega(\rho)\Big(U_tU_t^*AU_tU_t^*\Big), \end{eqnarray} for all $t>0$ and $A \in B(H)$, then
$\omega_t(I + \hat{\Lambda}\omega_t)^{-1}$ is a completely positive contraction
from $B(K)_*$ into $B(H)_*$ for every $t>0$.

Having seen that every $CP$-flow has an associated boundary weight map,
we naturally ask when a given map $\rho \rightarrow \omega(\rho)$
from $B(K)_*$ to $\mathfrak{A}(H)_*$ is the boundary weight
map of a $CP$-flow.  The answer is that if $\rho \rightarrow \omega(\rho)$ is a completely positive
map from
$B(K)_*$ into $\mathfrak{A}(H)_*$ satisfying
$\omega(\rho)(I-\Lambda(I_K)) \leq \rho(I_K)$ for all positive $\rho \in B(K)_*$, and
if $\omega_t(I + \hat{\Lambda} \omega_t)^{-1}$
is a completely positive contraction of $B(K)_*$ into $B(H)_*$ for every $t>0$,
then $\rho \rightarrow \omega(\rho)$ is the boundary weight map of a unique
$CP$-flow over $K$ (see Theorem 3.3 of \cite{bigpaper}).  This $CP$-flow is unital if and only if
$\omega(\rho)(I-\Lambda(I_K)) = \rho(I_K)$ for all $\rho \in
B(K)_*$.

Suppose $\alpha$ is a $CP$-flow over $\C$.  We identify its boundary
weight map with the single positive boundary weight
$\omega:=\omega(1) \in \mA$.  From above, $\omega$ has the form
$$\omega(\sqrt{I - \Lambda(1)} B \sqrt{I - \Lambda(1)}) = \sum_{i=1}^n (f_i, B f_i)$$
for some mutually orthogonal nonzero $L^2$-functions $\{f_k\}_{i=1}^n$
and unique $n \in \mathbb{N}\cup \{\infty\}$.  If $\alpha$ is unital, then
$\sum_{i=1}^n ||f_i||^2 = 1$,
and we say $\omega$ is \textit{normalized}.    We say $\omega$
is \textit{bounded} if there exists an $r>0$ such that
$|\omega(A)| \leq r ||A||$ for all $A \in \mathfrak{A}(L^2(0, \infty))$.  Otherwise, we say
$\omega$ is \textit{unbounded}.
From \cite{hugepaper}, we know that if $\omega$ is bounded, then the Bhat dilation $\alpha^d$ of
$\alpha$ is of type I$_n$, while if $\omega$ is unbounded, then $\alpha^d$ is of type
II$_0$.   Being type II$_0$ means that $\alpha_t^d$ is a proper $*$-endomorphism for all $t>0$
and that $\alpha^d$ has exactly one unit
$V=\{V_t\}_{t \geq 0}$ up to exponential scaling.  In other words, a semigroup
of bounded operators $W=\{W_t\}_{t \geq 0}$ acting on $H$ is a unit
for $\alpha^d$ if and only if, for some $\lambda \in \C$, we have
$W_t = e^{\lambda t}V_t$ for all $t \geq 0$.
This paragraph leads us to make the definition:
\begin{defn}
A normalized positive boundary weight $\nu \in \mA$ is said to be a
\emph{type I} (respectively, \emph{type II}) \emph{Powers weight} if
$\nu$ is bounded (respectively, unbounded).
\end{defn}

If $dim(K)>1$, we can naturally construct type II$_0$ $E_0$-semigroups by
combining type II Powers weights with $q$-positive maps acting on $M_n(\C)$ (Proposition
3.2 and Corollary 3.3 of \cite{Me}):

\begin{prop} \label{bdryweight} Let $H =\C^n \otimes L^2(0, \infty)$.
Let $\phi: M_n(\C) \rightarrow M_n(\C)$ be a unital completely
positive map with no negative eigenvalues, and let $\nu$ be a type
II Powers weight.  Let $\Omega_\nu: \mathfrak{A}(H) \rightarrow M_n(\C)$
be the map that sends $A=(A_{ij}) \in M_n(\mathfrak{A}(L^2(0, \infty))) \cong \mathfrak{A}(H)$
to the matrix $(\nu(A_{ij})) \in M_n(\C)$.

Then the map $\rho \rightarrow
\omega(\rho)$ from $M_n(\C)^*$ into
$\mathfrak{A}(H)_*$ defined by
$$\omega (\rho) (A) = \rho\Big(\phi(\Omega_\nu(A))\Big)$$
is the boundary weight map of a unital $CP$-flow $\alpha$ over $\C^n$ if and only if $\phi$ is
$q$-positive, in which case the Bhat minimal dilation $\alpha^d$ of $\alpha$ is a type
II$_0$ $E_0$-semigroup.
\end{prop}
In the notation of this proposition, we say $\alpha^d$ is the $E_0$-semigroup induced by the boundary
weight double $(\phi, \nu)$.  There is no ambiguity in doing so, since
$\alpha^d$ is unique up to conjugacy by Bhat's theorem.  Suppose that $(\phi, \nu)$ and $(\psi, \mu)$ are boundary weight doubles
which induce $E_0$-semigroups $\alpha^d$ and $\beta^d$.  When are $\alpha^d$
and $\beta^d$ cocycle conjugate?  We have a partial answer, and it involves
the following definition:

\begin{defn}
Let $\phi: M_n(\C) \rightarrow M_n(\C)$ and $\psi:
M_k(\C) \rightarrow M_k(\C)$ be $q$-positive maps.  We say a linear map $\gamma:
M_{n \times k}(\C) \rightarrow M_{n \times k}(\C)$
a \emph{corner} from $\phi$ to $\psi$ if the map
\begin{displaymath} \Upsilon \left( \begin{array}{cc} A_{n \times n} &
B_{n \times k} \\ C_{k \times n} & D_{k \times k}
\end{array} \right)  = \left( \begin{array}{cc} \phi(A_{n \times n}) &
\gamma (B_{n \times k})  \\
\gamma^* (C_{k \times n}) & \psi(D_{k \times k}) \\
\end{array} \right)
\end{displaymath}
is completely positive.  We say $\gamma$ is a \emph{$q$-corner} if
$\Upsilon$ is $q$-positive.  A $q$-corner $\gamma$ is called \emph{hyper maximal} if,
whenever
\begin{displaymath} \Upsilon \geq_q \Upsilon'  = \left( \begin{array}{cc} \phi' &  \gamma  \\
\gamma^* & \psi'\\
\end{array} \right) \geq_q 0,
\end{displaymath}
we have $\Upsilon = \Upsilon'$.
\end{defn}

The main result of \cite{Me} with regard to comparing
$E_0$-semigroups induced by boundary weight doubles $(\phi, \nu)$
and $(\psi, \nu)$ is the following, which unfortunately requires
$\nu$ to have a very specific form:

\begin{prop}\label{hypqc}
Let $\phi: M_n(\C) \rightarrow M_n(\C)$ and $\psi: M_k(\C)
\rightarrow M_k(\C)$ be unital $q$-positive maps, and let $\nu$ be a
type II Powers weight of the form $$\nu(\sqrt{I -
\Lambda(1)}B\sqrt{I - \Lambda(1)})=(f, Bf).$$ The boundary weight
doubles $(\phi, \nu)$ and $(\psi, \nu)$ induce cocycle conjugate
$E_0$-semigroups if and only if there is a hyper maximal $q$-corner
from $\phi$ to $\psi$.
\end{prop}

Let $\phi: M_n(\C) \rightarrow M_n(\C)$ be unital and $q$-positive, and let
$U \in M_n(\C)$ be unitary.  Define $\phi_U: M_n(\C) \rightarrow M_n(\C)$ by
$$\phi_U(A)=U^*\phi(UAU^*)U$$ for all $A \in M_n(\C)$.  It is straightforward
to show that $\phi_U$ is also unital and $q$-positive.  We note that the map $\gamma: M_n(\C)
\rightarrow M_n(\C)$ defined by $\gamma(A)=\phi(AU^*)U$ is a hyper maximal
$q$-corner from $\phi$ to $\phi_U$.  Indeed, it is easy to check that $\gamma$
is a $q$-corner from $\phi$ to $\phi_U$ (see Proposition 4.5 of \cite{Me}).
To see that $\gamma$ is hyper maximal,
we observe that if
\begin{displaymath}
\left( \begin{array}{cc}
\phi & \gamma \\ \gamma^* & \phi_U
\end{array} \right) \geq_q
\left( \begin{array}{cc}
\phi' & \gamma \\ \gamma^* & \phi'_U
\end{array} \right) \geq_q 0,
\end{displaymath}
then $\phi'(I) \leq I$ and $\phi'_U(I) \leq I$, yet
\begin{displaymath}
\left( \begin{array}{cc}
\phi'(I) & \gamma(U) \\ \gamma^*(U^*) & \phi'_U(I)
\end{array} \right) =
\left( \begin{array}{cc}
\phi'(I) & U \\ U^* & \phi'_U(I)
\end{array} \right) \geq 0,
\end{displaymath}
hence $\phi'(I)=\phi'_U(I)=I$.  But $\phi - \phi'$ and $\phi_U - \phi'_U$
are completely positive, so $$||\phi-\phi'||=||\phi(I)-\phi'(I)||=0=||\phi_U(I)
-\phi'_U(I)=||\phi_U - \phi'_U||,$$ thus $\phi=\phi'$ and $\phi_U= \phi_U'$.  This
shows that $\gamma$ is hyper maximal, whereby Proposition \ref{hypqc} gives us the following:

\begin{prop}\label{dumb}
Let $\phi: M_n(\C) \rightarrow M_n(\C)$ be unital and $q$-positive,
and let $U \in M_n(\C)$ be unitary.  If $\nu$ is a type II Powers
weight of the form $$\nu(\sqrt{I - \Lambda(1)}B\sqrt{I -
\Lambda(1)})=(f, Bf),$$ then $(\phi, \nu)$ and $(\phi_U, \nu)$
induce cocycle conjugate $E_0$-semigroups.
\end{prop}

In fact, if $\nu$ is an arbitrary type II Powers weight,
then $(\phi, \nu)$ and $(\phi_U, \nu)$ induce \textit{conjugate}
$E_0$-semigroups (\cite{DJ}). We will not use this result here,
except as justification for the following definition.
%

\begin{defn}  Let $\phi: M_n(\C) \rightarrow M_n(\C)$ be $q$-positive.
We say $\psi$ is \emph{conjugate} to $\phi$ if $\psi=\phi_U$ for some unitary
$U \in M_n(\C)$.
\end{defn}

In other words, $\psi$ is conjugate to $\phi$ if and only if
there is a $*$-isomorphism $\theta:
M_n(\C) \rightarrow M_n(\C)$ such that $\theta \circ \psi \circ \theta^{-1} = \phi.$  This
is analogous to the notion of conjugacy for $E_0$-semigroups and is appropriate
in light of the preceding paragraph.
We recall the classification of all $q$-pure maps $\phi: M_n(\C)
\rightarrow M_n(\C)$ which are rank one or invertible, along
with the main cocycle conjugacy results of \cite{Me}
(see Lemma 5.2 and Theorem 5.4 of \cite{Me}):
\begin{thm}\label{statesbig}  A unital rank one linear map $\phi: M_n(\C)
\rightarrow M_n(\C)$ is $q$-positive if and only if it has the form
$\phi(A)=\rho(A)I$ for some state $\rho \in M_n(\C)^*$.  Such a map $\phi$
is $q$-pure if and only if $\rho$ is faithful.

Let $\phi$ and $\psi$ be unital rank one $q$-pure maps on
$M_n(\C)$ and $M_k(\C)$, respectively, and let $\nu$ be a type II Powers weight
of the form
$\nu(\sqrt{I-\Lambda(1)}B\sqrt{I-\Lambda(1)})=(f,Bf)$.  The boundary weight doubles $(\phi, \nu)$ and
$(\psi, \nu)$ induce cocycle conjugate $E_0$-semigroups if and
only if $n=k$ and $\phi$ is conjugate to $\psi$.
\end{thm}
Furthermore, if $\nu$ and $\mu$ are type II Powers weights and
$\phi$ and $\psi$ are rank one unital $q$-pure maps on $M_n(\C)$
and $M_k(\C)$, respectively, then $(\phi, \nu)$ and $(\psi, \mu)$
cannot induce cocycle conjugate $E_0$-semigroups unless there is a
corner $\gamma$ from $\phi$ to $\psi$ such that $||\gamma||=1$
(Lemma 5.3 of \cite{Me}). A consequence of this result is that if
$n>1$, then none of the $E_0$-semigroups induced by boundary weight
doubles $(\phi, \nu)$ for unital rank one $q$-pure $\phi: M_n(\C)
\rightarrow M_n(\C)$ and $\nu$ 
of the form $\nu(\sqrt{I-\Lambda(1)}B\sqrt{I-\Lambda(1)})=(f,Bf)$ are cocycle conjugate to any of the
$E_0$-semigroups constructed by Powers in the case that $dim(K)=1$
in \cite{bigpaper}. However, for $q$-pure maps that are invertible
rather than rank one, the opposite holds (Theorems 6.11 and 6.12 of
\cite{Me}):

\begin{thm}\label{inverts}  An invertible unital linear map $\phi: M_n(\C) \rightarrow M_n(\C)$
is $q$-positive if and only if $\phi^{-1}$ is conditionally negative, and $\phi$ is $q$-pure
if and only if
$\phi^{-1}$ is of
the form
$$\phi^{-1}(A) = A + YA + AY^*$$ for some $Y \in M_n(\C)$ with $Y=-Y^*$ and
$tr(Y)=0$. Equivalently, $\phi$ is $q$-pure if and only if it is
conjugate to a Schur map $\psi$ that satisfies

\begin{equation*} \psi(a_{jk}e_{jk}) = \left(
\begin{array}{cc}
\frac{a_{jk}}{1+i(\lambda_j - \lambda_k)}e_{jk} & \textrm{if } j<k \\
a_{jk}e_{jk} & \textrm{if } j=k \\
\frac{a_{jk}}{1-i(\lambda_j - \lambda_k)}e_{jk}& \textrm{if } j>k
\end{array} \right)
\end{equation*}
for all $j,k=1,\ldots, n$ and all $A=\sum a_{ij}e_{ij} \in M_n(\C)$,
where $\lambda_1, \ldots, \lambda_n \in \R$ and $\sum_{j=1}^n
\lambda_j = 0$.

If $\nu$ is a type II Powers weight of the form
$\nu(\sqrt{I-\Lambda(1)}B\sqrt{I-\Lambda(1)})=(f,Bf)$, then the $E_0$-semigroup
induced by $(\phi,\nu)$ is cocycle conjugate to the $E_0$-semigroup induced by $(\imath_\C, \nu)$
for $\imath_\C$ the identity map on $\C$ (this is the $E_0$-semigroup induced by $\nu$ in the sense of
\cite{bigpaper}).
\end{thm}

\section{$\mathcal{E}_n$ and the limiting map $L_\phi$}
Suppose $\phi: M_n(\C) \rightarrow M_n(\C)$ is a $q$-positive map and $||t
\phi(I + t \phi)^{-1}|| < 1$ for all $t>0$. In \cite{Me}, we saw that we could form a limit
$L_\phi= \lim_{t \rightarrow \infty} t \phi(I + t \phi)^{-1}$. This
limiting process was the key to classifying the rank one $q$-pure
maps acting on $M_n(\C)$.  We begin this section by revisiting $L_\phi$:

\begin{lem}\label{limit}  Suppose $\phi: M_n(\C) \rightarrow M_n(\C)$ is a
non-zero $q$-positive map such that $||t \phi(I + t \phi)^{-1}||<1$ for all $t
\geq 0$. Then the maps $t \phi(I + t \phi)^{-1}$ have a unique limit
$L_\phi$ as $t \rightarrow \infty$, and $||L_\phi||=1$. Furthermore, $L_\phi$
is completely positive, $L_\phi \circ \phi=
\phi \circ L_\phi= \phi$, $\range(L_\phi) =
\range(\phi)$, $\nullspace(\phi)=\nullspace(L_\phi)$, and $L_\phi ^2 = L_\phi$.
\end{lem}
{\bf Proof:} A compactness argument shows that since $||t \phi(I + t \phi)^{-1}||<1$
for all $t>0$, the maps $t \phi(I + t \phi)^{-1}$
have some norm limit $L_\phi$ as $t \rightarrow \infty$, where $||L_\phi|| \leq 1$.  To see this limit is
unique, we let $M \in M_{2n}(\C)$ be the matrix for $\phi$ with respect to some
orthonormal basis of $M_n(\C)$ and note that the entries of $tM (I +
t M)^{-1}$ are (necessarily bounded) rational functions of $t$ and
thus each have unique limits.  $L_\phi$ is completely positive since it
is the norm limit of completely positive maps.

For every $t>0$, let $M_t = (I + t \phi)/t$, so $M_t \rightarrow
\phi$ as $t \rightarrow \infty$. Given any $A \in M_n(\C)$, we find
\begin{eqnarray*} \phi(L_\phi(A)) & = & \lim_{t \rightarrow \infty} M_t
(t \phi(I + t \phi)^{-1})(A)) = \lim_{t \rightarrow \infty}
\Big(\frac{I + t \phi}{t}\Big) t\phi(I + t \phi)^{-1}(A) \\ & = &
\phi(I + t \phi)(I + t \phi)^{-1}(A) = \phi(A),
\end{eqnarray*}
so $\phi \circ L_\phi = \phi$.  But $L_\phi$ commutes with $\phi$ by
construction, so $L_\phi \circ \phi = \phi$, hence $\range(\phi)
\subseteq \range(L_\phi)$. Trivially $\range(L_\phi) \subseteq
\range(\phi)$, so $\range(\phi)= \range(L_\phi)$, whereby the fact that $L_\phi \circ
\phi = \phi$ implies that $L_\phi$ fixes its range, hence
$L_\phi^2 = L_\phi$.  Since $L_\phi$ is a nonzero contraction,
we have $||L_\phi||=1$.  To finish the proof, we need only show
that $\phi$ and $L_\phi$ have the same nullspace.  The fact that
$\nullspace(L_\phi)\subseteq \nullspace(\phi)$ follows trivially from the established
equality $\phi= \phi \circ L_\phi$.  On the other hand,
if $\phi(A)=0$, then $(I + t\phi)^{-1}(A)=A$ for all $t \geq 0$, hence
$$L_\phi(A) = \lim_{t \rightarrow \infty} t \phi(I +
t \phi)^{-1}(A)=\lim_{t \rightarrow \infty} t \phi(A)=0.$$ \qqed

Any unital $q$-positive $\phi: M_n(\C) \rightarrow M_n(\C)$
satisfies the conditions of the above lemma, since for all $t \geq
0$ we have $$||t \phi(I+t\phi)^{-1}||= ||t \phi(I+t\phi)^{-1}(I)|| =
\frac{t}{1+t}.$$

\begin{defn} For each
$n \in \mathbb{N}$, let $\mathcal{E}_n$ be
the set of all unital completely positive maps $\Phi: M_n(\C)
\rightarrow M_n(\C)$ such that $\Phi^2=\Phi$.
\end{defn}
If $\phi: M_n(\C) \rightarrow M_n(\C)$ is unital and $q$-positive,
then $L_\phi \in \mathcal{E}_n$ by Lemma \ref{limit}.  On the other hand,
let $\Phi \in \mathcal{E}_n$ be arbitrary.  Since
$\Phi^2=\Phi$,  it follows that $I + t \Phi$ is invertible for all $t \geq 0$
and $t \Phi(I + t\Phi)^{-1}= (t/(1+t)) \Phi$, so $\Phi$ is $q$-positive
and $\Phi= L_\Phi$.  Therefore,
$$\mathcal{E}_n = \{ L_\phi \ | \  \phi: M_n(\C) \rightarrow M_n(\C), \
\phi(I)=I,  \textrm{ and } \phi \geq_q 0 \}.$$  Note that membership
in $\mathcal{E}_n$ is invariant under conjugacy:  If $\Phi \in
\mathcal{E}_n$ and $U \in M_n(\C)$ is unitary, then $\Phi_U$ is unital and
completely positive by construction, and $\Phi_U^2=\Phi_U$ since
\begin{eqnarray*} \Phi_U ^2(A) & = & \Phi_U (U^* \Phi(UAU^*)U)=U^*\Phi \Big[ U
\Big(U^* \Phi(UAU^*)U \Big) U^*\Big]U \\ &=& U^* \Phi^2(UAU^*)U =
U^* \Phi(UAU^*)U=\Phi_U(A) \end{eqnarray*} for all $A \in M_n(\C)$.

The rest of this section is devoted to classifying the elements of
$\mathcal{E}_2$ and $\mathcal{E}_3$ up to conjugacy. As we will see
in the next section, this is a key step in classifying all unital
$q$-positive maps $\phi: M_2(\C) \rightarrow M_2(\C)$ and in showing
that a large class of $q$-positive maps acting on $M_3(\C)$ cannot
be $q$-pure.
\\

{\bf Remark:}  It is possible for a unital completely positive map $\phi: M_2(\C) \rightarrow M_2(\C)$
to have rank 3. For example, define $\phi: M_2(\C) \rightarrow M_2(\C)$ by
\begin{displaymath}
\phi(A)= \frac{1}{3} \left( \begin{array}{cc} 2 a_{11} + a_{22} & a_{12}+a_{21} \\
a_{12}+a_{21} & a_{11}+2a_{22} \end{array} \right).
\end{displaymath}
We see that $\phi$ has rank $3$, since
\begin{displaymath}
\range(\phi)= \Big\{ \left( \begin{array}{cc} a & b \\
b & c \end{array} \right):  a, b, c \in \C \Big\}.
\end{displaymath}
Furthermore, $\phi$ is completely positive since it is the sum of completely positive maps,
as
$$\phi(A)= \frac{1}{3} \Big( A + SAS^* + D(A)\Big),$$
where $S= e_{12} + e_{21}$ and $D$
is the diagonal map $D(A)=a_{11}e_{11} + a_{22}e_{22}$.  However, it turns out that $\phi$ is not $q$-positive.
In fact, we will see from our classification of $\mathcal{E}_2$ that
no unital $q$-positive map $\phi$
acting on $M_2(\C)$ can have rank 3:

\begin{prop}\label{E2}
Let $\Phi: M_2(\C) \rightarrow M_2(\C)$ be a unital linear map. Then $\Phi \in \mathcal{E}_2$
if and only if, up to conjugacy, $\Phi$ has one of the forms below:
\begin{enumerate}[(i)]
\item $\Phi(A)=\rho(A)I$ for all $A \in M_2(\C)$, where $\rho \in M_2(\C)^*$ is a state.
\item $\Phi(A)= a_{11}e_{11} + a_{22}e_{22}$ for all $A=\sum_{i,j=1}^2 a_{ij}e_{ij} \in M_2(\C)$.
\item $\Phi(A)=A$ for all $A \in M_2(C)$.
\end{enumerate}
Consequently, if $\phi: M_2(\C) \rightarrow M_2(\C)$ is unital and $q$-positive,
then $\rank(\phi) \neq 3$.
\end{prop}
{\bf Proof:}  By inspection, maps $(i)$ through $(iii)$ (and therefore
their conjugates) are in $\mathcal{E}_2$.
On the other hand, suppose $\Phi$ is an element of $\mathcal{E}_2$.  If
$\Phi$ has rank one, then it trivially has the form $(i)$.  If $\rank(\Phi)
\geq 2$, then $\Phi(I)=I$ and
$\Phi(A_1)=A_1$ for some $A_1$ linearly independent from $I$.
Since $\Phi$ is completely positive and thus self-adjoint in the
sense that $\Phi(A^*)=\Phi(A)^*$ for all $A$, we have
$\Phi(A_1+A_1^*)=A_1+A_1^*$ and
$\Phi(i(A_1-A_1^*))=i(A_1-A_1^*)$.  A quick exercise in linear algebra
shows that
the self-adjoint matrices $A_1+A_1^*$ and $i(A_1-A_1^*)$ cannot
both be multiples of $I$, whereby we conclude that $\Phi(M)=M$ for some self-adjoint
$M \in M_n(\C)$ linearly independent from $I$.

Letting $U$ be a unitary matrix such that $U^*MU=D$ for some
diagonal matrix $D$, we note that $D$ is linearly independent from
$I$. We observe that
$\Phi_U(I)=I$ and $\Phi_U(D) = U^*\Phi(UDU^*)U=U^*MU=D$, which implies
$\Phi_U(e_{11})=e_{{11}}$ and
$\Phi_U(e_{22})=e_{22}$.  We claim that $\Phi_U (e_{12})= b
e_{12}$ for some $b \in \C$. Indeed, write

\begin{displaymath}\Phi_U (e_{12})= \left( \begin{array}{cc} a
 & b \\
c & d
\end{array} \right).
\end{displaymath}
Since $\Phi_U$ is $2$-positive, we have
\begin{displaymath}0 \leq \left( \begin{array}{cc} \Phi_U(e_{11}) &
\Phi_U(e_{12})  \\ \Phi_U(e_{21}) & \Phi_U(e_{22})
\end{array} \right) =
\left( \begin{array}{cccc} 1 & 0 & a & b  \\
0 & 0 & c & d  \\ \overline{a} & \overline{c} & 0 & 0 \\
\overline{b} & \overline{d} & 0 & 1
\end{array} \right).
\end{displaymath}
Positivity of the above matrix implies $a=c=d=0$, hence
$\Phi_U(e_{12})=be_{12}$ and
$\Phi_U(e_{21})=\overline{b}e_{21}$.  Therefore $\Phi_U$ is
merely the Schur mapping
\begin{displaymath}\left( \begin{array}{cc}
a_{11}
 & a_{12} \\
a_{21} & a_{22} \end{array} \right) \rightarrow \left(
\begin{array}{cc} a_{11}
 & b a_{12} \\
\overline{b} a_{21} & a_{22} \end{array} \right). \end{displaymath}
Since $(\Phi_U)^2=\Phi_U$ we have $b^2=b$, so either $b=0$ (in which case $\Phi_U$
has the form $(ii)$) or $b=1$ (in which case $\Phi_U$ is the identity map $(iii)$).

For the final statement of the theorem, we note that if $\phi:
M_2(\C) \rightarrow M_2(\C)$ is unital and $q$-positive, then
$L_\phi \in \mathcal{E}_2$ and $\rank(\phi)=\rank(L_\phi)$, so
$\rank(L_\phi) \in \{1,2,4\}$ by what we have just shown. \qqed

We turn our attention to classifying the elements of $\mathcal{E}_3$
up to conjugacy. Our task is made much easier by the fact that each
of its elements with rank greater than one must destroy or fix a
rank one projection:

\begin{lem}\label{fixordestroy} Suppose $\Phi \in \mathcal{E}_3$ and $\rank(\Phi)>1$.  If $\Phi$ does not
annihilate any nonzero projections, then $\Phi$ fixes some rank one
projection $E$.
\end{lem}

{\bf Proof:}  Since $rank(\Phi) \geq 2$ and $\Phi$ fixes its range,
$\Phi$ fixes some $M \in M_3(\C)$ linearly independent from $I$.
Arguing as in the proof of Proposition \ref{E2}, we may assume $M=M^*$,
and of course we may assume $||M|| = 1$.  Since $M$ is self-adjoint and has norm
one, we know that at least one of the numbers $1$ and $-1$ is an eigvenvalue
of $M$.  Therefore, replacing $M$ with $-M$ if necessary, we may assume that
$1$ is an eigenvalue of $M$.  Diagonalizing $M$ by a unitary $U \in M_3(\C)$
so that the eigenvalues of $D:=UMU^*$ are listed in decreasing order,
we have
\begin{displaymath}
D= \left( \begin{array}{ccc} 1 & 0 & 0 \\ 0 & \lambda_1 & 0
\\ 0 & 0 & \lambda_2 \end{array} \right)
\end{displaymath}
where $1 \geq \lambda_1 \geq \lambda_2$.  Note that $\lambda_2 \neq 1$ since $D \neq I$.
Since $\Phi_U$ fixes (respectively, annihilates) a projection $P$ if and only if $\Phi$ fixes
(respectively, annihilates) the projection $UPU^*$, it suffices to show that $\Phi_U$ fixes
a rank one projection.

Note that $\Phi_U(I)=I$ and $\Phi_U(D)=D$, so

\begin{equation} \label{mon}
\Phi_U(I-D)= I- D= \left( \begin{array}{ccc} 0 & 0 & 0 \\ 0 & 1-\lambda_1 & 0
\\ 0 & 0 & 1-\lambda_2 \end{array} \right) \geq 0.
\end{equation}
If $\lambda_1=1$, then $\Phi_U$ fixes $e_{33}$ and the lemma follows.  If $\lambda_1 \neq 1$,
then we let
$b= (1-\lambda_2)/(1-\lambda_1)>0$.  By complete positivity
of $\Phi_U$ and equation \eqref{mon},
\begin{equation} \label{Fix1}
0 \leq \Phi_U(e_{22}) \leq \Phi_U(e_{22} + b e_{33}) = e_{22} + b e_{33}.\end{equation}
We also note that
\begin{displaymath}
\Phi_U(D-\lambda_2 I)= D- \lambda_2I = \left( \begin{array}{ccc} 1-\lambda_2 & 0
& 0 \\ 0 & \lambda_1-\lambda_2  & 0
\\ 0 & 0 & 0 \end{array} \right) \geq 0.
\end{displaymath}
If $\lambda_1=\lambda_2$, then $\Phi_U$ fixes $e_{11}$ and the lemma follows.  If
$\lambda_1 \neq \lambda_2$, we let $c=(1-\lambda_2)/(\lambda_1-\lambda_2)>0$ and note that
\begin{equation} \label{Fix2}
0 \leq \Phi_U(e_{22}) \leq \Phi_U(c e_{11} + e_{22}) = ce_{11} + e_{22}.\end{equation}
Equation \eqref{Fix1} implies that the $11$ entry of $\Phi_U(e_{22})$ is zero,
while equation \eqref{Fix2} implies that the $33$ entry of $\Phi_U(e_{22})$ is zero.
Therefore, $\Phi_U(e_{22})= \lambda e_{22}$ for some $\lambda \geq 0$.  Since $\Phi_U^2
=\Phi_U$ we have $\lambda \in \{0,1\}$, whereby the fact that $\Phi_U$
does not annihilate any nonzero projections implies $\lambda=1$.  Thus,
$\Phi_U$ fixes $e_{22}$, so $\Phi$ fixes the rank one projection $U e_{22} U^*$.

\qqed

Before proceeding further, we will need the following two standard results regarding
completely positive maps:
\begin{lem}\label{FAF} Let $K$ be a separable Hilbert space, and let
$\phi: B(K) \rightarrow B(K)$ be a normal completely positive map.
If $\phi(E)=0$ for some projection $E$, then $\phi(A)=\phi(FAF)$ for
all $A \in B(K)$, where $F=I-E$.
\end{lem}
From \cite{arveson}, we know that $\phi$ can be written in the form
$\phi(A)=\sum_{i=1}^p S_i A S_i^*$ for some $p \in \mathbb{N} \cup
\{\infty\}$ and operators $\{S_i\}_{i=1}^p$ in $B(K)$.  We note that
$$0= \phi(E)= \sum_{i=1}^p S_i E S_i^* = \sum_{i=1}^p S_i EES_i^*=
\sum_{i=1}^p (S_i E)(S_iE)^*,$$ so
$S_iE=0=ES_i^*$ for all $i$.  Therefore, $\phi(EAE)=
\phi(EAF)=\phi(FAE)=0$ for all $A \in B(K)$, hence
\begin{eqnarray} \nonumber \phi(A) &=& \phi((E+ F)A(E+F)) \\
& = & \nonumber \phi(EAE)+ \phi(EAF)+\phi(FAE)+\phi(FAF)
\\ \nonumber & = & \phi(FAF).  \end{eqnarray} \qqed

\begin{lem}\label{form}
Let $K$ be a separable Hilbert space, and let
$\phi: B(K) \rightarrow B(K)$ be a normal and unital
completely positive map.  Suppose $\phi$
fixes a projection $E$.  Then
$$\phi(A)= E\phi(EAE)E + E\phi(EAF)F + F \phi(FAE)E + F\phi(FAF)F$$
for all $A$, where $F=I-E$.
\end{lem}
{\bf Proof:}  By hypothesis, we can write $\phi$ in the form
$\phi(A)= \sum_{i=1}^p S_iAS_i^*$ .  Since $\phi(I)=I$ and
$\phi(E)=E$, we have $\phi(F)=\phi(I-E)=I-E=F$.  Therefore, $S_i
ES_i^* \leq E$ and $S_iFS_i^* \leq F$ for all $i$.  Note that $E S_i
F= FS_i^*E = 0$ for all $i$, since
$$(E S_iF)(E S_i F)^* = E S_i F F S_i^* E = E(S_i F S_i^*)E \leq EFE = 0.$$  An analogous
argument shows that $FS_iE=  ES_i^*F=0$ for all $i$.  Writing
$$\phi(A) = (E+F)\phi\Big((E+F)A(E+F)\Big) (E+F)$$ and expanding the
right hand side using the above makes most of the terms vanish,
yielding the result. \qqed

With the previous three lemmas in hand, we are able to classify the
elements of $\mathcal{E}_3$ in two steps.  

\begin{lem}\label{annih2} Let $\Phi: M_3(\C) \rightarrow M_3(\C)$ be a
unital map such that $\Phi(E)=0$ for some nonzero projection $E$.
Then $\Phi \in \mathcal{E}_3$ if and only if, up to conjugacy,
$\Phi$ has one of the following forms for some $\lambda \in [0,1]$:
\begin{equation} \label{arg1} \tag{I} \Phi(A) =
\Big(\lambda a_{22} + (1-\lambda)a_{33}\Big)I;
\end{equation}

\begin{equation} \label{arg3} \tag{II} \Phi(A) = \left( \begin{array}{ccc}
\lambda a_{22} + (1-\lambda)a_{33} & 0 & 0
\\ 0 & a_{22} & 0 \\ 0 & 0 & a_{33}
\end{array} \right);
\end{equation}

\begin{equation} \label{arg5} \tag{III} \Phi(A) = \left( \begin{array}{ccc}
\lambda a_{22} + (1-\lambda)a_{33} & 0 & 0
\\ 0 & a_{22} & a_{23} \\ 0 & a_{32} & a_{33}
\end{array} \right);
\end{equation}
\end{lem}
{\bf Proof:}  The backward direction follows from inspection of the
maps \eqref{arg1} through \eqref{arg5}. For the forward direction,
suppose $\Phi \in \mathcal{E}_3$ and $\Phi(E)=0$.  Let $E'$ be any
rank one subprojection of $E$, observing that $\Phi(E')=0$.
Unitarily diagonalizing $E'$ so that $ Z^* E' Z = e_{11}$, we have
$\Phi_Z (e_{11})=0$.  By Lemma \ref{FAF} it follows that
$\Phi_Z(A)=\Phi_Z(FAF)$ for all $A \in M_3(\C)$, where $F=I-e_{11}$.
Replacing $\Phi_Z$ with $\Phi$ (as we are only concerned with $\Phi$
up to conjugacy), we write
\begin{displaymath}
\Phi (A) = \left(   \begin{array}{ccc} \tau_1(A) & \tau_2(A) &
\tau_3(A)
\\ \tau_2^*(A) & \multicolumn{2}{c}{\multirow{2}{*}{$[\Psi(A)]$}}
\\ \tau_3^*(A) \end{array} \right)
\end{displaymath}
for some linear functionals $\tau_j$ ($j=1, 2, 3$) and some map
$\Psi: M_3(\C) \rightarrow M_2(\C)$.  But $\Phi(A)=\Phi(FAF)$ for all $A \in
M_3(\C)$, so $\Phi(e_{1j})=\Phi(e_{j1})=0$ for $j=1,2,3$.
Therefore, for every $A \in M_3(\C)$, $\Psi(A)$ and each $\tau_j(A)$ depend only on
the bottom right $2 \times 2$ minor of $A$.  In other words, if we let
\begin{displaymath}
G= \left( \begin{array}{ccc} 0 & 1 & 0 \\ 0 & 0 & 1 \end{array}
\right) \in M_{2 \times 3}(\C)
\end{displaymath}
and define $\rho_j \in M_2(\C)^*$ ($j=1,2,3$) and $\psi: M_2(\C) \rightarrow M_2(\C)$
by $\rho_j(B)=\tau(G^*BG)$ and $\psi(B)=\Psi(G^*BG)$, then for all $A \in M_3(\C)$,
\begin{equation}\label{G} \Phi(A)= \left(\begin{array}{ccc} \rho_1(GAG^*) &
\rho_2(GAG^*) & \rho_3(GAG^*) \\ \rho_2 ^*(GAG^*) &
\multicolumn{2}{c}{\multirow{2}{*}{$[\psi(GAG^*)]$}}
\\ \rho_3^*(GAG^*)
\end{array} \right).
\end{equation}
Note that $\psi(B) = G \Phi(G^*BG)G^*$ for all $B \in M_2(\C)$, so
$\psi$ is completely positive, and $\psi$ is unital since for the
identity matrix $I_2 \in M_2(\C)$, we have
$$\psi(I_2)=G \Phi(G^*G)G^* = G \Phi(F) G^*= G I G^*= I_2.$$  Furthermore, $\psi^2=\psi$, since
\begin{eqnarray*}
\psi^2(B) & = & \psi\Big(G \Phi(G^*BG)G^*\Big)= G \Phi\Big(G^*G \Phi(G^*BG)G^*G
\Big)G^* \\ & = &
G \Phi\Big(F \Phi(G^*BG) F\Big) G^* = G \Phi\Big(\Phi(G^*BG)
\Big)G^* \\ & = & G \Phi(G^*BG)G^* = \psi(B),
\end{eqnarray*}
where for the fourth equality we used the fact that
$\Phi(A)=\Phi(FAF)$ for all $A \in M_3(\C)$.  Therefore, $\psi \in
\mathcal{E}_2$, whereby Proposition \ref{E2} implies that
$rank(\psi) \neq 3$.

{\bf Case (i):} If $rank(\psi)=1$, then $\psi$ is of the form
$\psi(B)=\rho(B)I_2$, where $\rho \in M_2(\C)^*$ satisfies
$\rho(I_2)=1$. By equation \eqref{G} and the fact that
$\Phi^2=\Phi$, we have
\begin{equation}\label{rhoj1} \rho_j(GAG^*)=\rho_j
(\psi(GAG^*))=\rho_j(\rho(GAG^*)I_2) = \rho(GAG^*)\rho_j(I_2)
\end{equation} for every $A \in M_3(\C)$ and $j=1,2,3$. But
$\Phi(I)=I$, so $\rho_1(I_2)=1$ while $\rho_2(I_2)=\rho_3(I_2)=0$,
so equation \eqref{rhoj1} implies $\rho_1= \rho$ and $\rho_2=\rho_3
\equiv 0$.

Since $\rho \in M_2(\C)^*$ is a state, there is some
$\lambda \in [0,1]$ and a
unitary matrix $S \in M_2(\C)$ such that
$\rho(SBS^*)=\lambda b_{11} + (1-\lambda)b_{22}$ for all $B \in M_2(\C)$.
Therefore,
\begin{equation*}\rho \Big(S \Big[GAG^*\Big] S^* \Big) = \lambda a_{22} + (1 -\lambda) a_{33}
\end{equation*}
for all $A \in M_3(\C)$.  Letting
\begin{displaymath}
R = \left( \begin{array}{ccc}
1 & 0 & 0 \\ 0 &
\multicolumn{2}{c}{\multirow{2}{*}{$[S]$}}
\\ 0
\end{array} \right),
\end{displaymath}
we see that $\Phi_R$ has the form \eqref{arg1}.

{\bf Case (ii):} If $rank(\psi)=2$, then Lemma \ref{E2} implies that
for some $2 \times 2$ unitary $V$, $\psi_V$ is the diagonal map
\begin{displaymath}\psi_V\left(\begin{array}{cc} b_{11} & b_{12}
\\ b_{21} & b_{22} \end{array} \right)
= \left(\begin{array}{cc} b_{11} & 0
\\ 0 & b_{22} \end{array} \right).
\end{displaymath}
Let $U \in M_3(\C)$ be the $3 \times 3$ unitary matrix
\begin{displaymath}
U = \left( \begin{array}{ccc}
1 & 0 & 0 \\ 0 &
\multicolumn{2}{c}{\multirow{2}{*}{$[V]$}}
\\ 0
\end{array} \right).
\end{displaymath}

Then $\Phi_U(e_{11}) = U^*\Phi(Ue_{11}U^*)U=U^*\Phi(e_{11})U=0$ and
$G \Phi_U(G^* BG)G^*=\psi_V(B)$ for all $B \in M_2(\C)$.  Therefore,
$\Phi_U$ has the form below for some linear functionals $\rho_j'$,
$j=1,2,3$:

\begin{eqnarray} \Phi_U(A)= \left(\begin{array}{ccc}
\rho_1'(GAG^*) & \rho_2'(GAG^*) & \rho_3'(GAG^*)
\\ \rho_2'^*(GAG^*) & a_{22} & 0 \\ \rho_3'^*(GAG^*) & 0 & a_{33} \end{array}
\right). \end{eqnarray}

Replacing $\Phi_U$ with $\Phi$ and erasing the primes on the functionals for simplicity
of notation, we have

\begin{eqnarray*} \Phi(A)= \left(\begin{array}{ccc}
\rho_1(GAG^*) & \rho_2(GAG^*) & \rho_3(GAG^*)
\\ \rho_2^*(GAG^*) & a_{22} & 0 \\ \rho_3^*(GAG^*) & 0 & a_{33} \end{array}
\right). \end{eqnarray*}

Positivity of the matrices $\Phi(e_{22})$ and $\Phi(e_{33})$ yields
\begin{equation} \label{h1}
\rho_3 \left( \begin{array}{cc} 1 & 0 \\ 0 & 0  \end{array} \right)=0 \ \textrm{ and }
\
\rho_2 \left( \begin{array}{cc} 0 & 0 \\ 0 & 1 \end{array} \right)=0,
\end{equation}
respectively.  Since $\Phi$ is unital we have
\begin{equation} \label{h2}
\rho_3 \left( \begin{array}{cc} 1 & 0 \\ 0 & 1  \end{array} \right)
 = \rho_2 \left( \begin{array}{cc} 1 & 0 \\ 0 & 1 \end{array} \right)=0,
\end{equation}
and combining equations \eqref{h1} and \eqref{h2} gives us
\begin{equation} \label{j1} \rho_2(D)=\rho_3(D)=0 \ \textrm { for all diagonal matrices }
D \in M_2(\C).
\end{equation}

  For $j=1,2,3$, the fact that $\Phi^2(e_{23})=\Phi(e_{23})$
implies
\begin{equation} \label{j2}
\rho_j \left( \begin{array}{cc} 0 & 1 \\ 0 & 0  \end{array} \right)
= \rho_j \left( \psi \Big[\begin{array}{cc} 0 & 1 \\ 0 & 0  \end{array}\Big]
 \right) = \rho_j \left( \begin{array}{cc} 0 & 0 \\ 0 & 0 \end{array} \right) =0,
\end{equation}
and similarly, since $\Phi^2(e_{32})=\Phi(e_{32})$, we have
\begin{equation} \label{j3}
\rho_j \left( \begin{array}{cc} 0 & 0 \\ 1 & 0  \end{array} \right)=0.
\end{equation}
From equations \eqref{j1}, \eqref{j2}, and \eqref{j3}, we have $\rho_2=\rho_3\equiv 0$
and
\begin{equation} \label{r1}
\rho_1 \left( \begin{array}{cc} b_{11} & b_{12} \\ b_{21} & b_{22}  \end{array} \right)
= \rho_1 \left( \begin{array}{cc} b_{11} & 0 \\ 0 & b_{22}  \end{array} \right)
\end{equation}
for all $B \in M_2(\C)$.  From equation \eqref{r1} and the fact that
$\Phi$ is unital, there is some $\lambda \in [0,1]$ such that
$$\rho_1(GAG^*)=\lambda a_{22} + (1-\lambda)a_{33}$$ for all $A \in M_3(\C)$,
hence $\Phi$ satisfies \eqref{arg3}.

{\bf Case (iii):} If $rank(\psi)=4$, then $\psi$ is the identity map
by Lemma \ref{E2}, so

\begin{eqnarray*} \Phi(A)= \left(\begin{array}{ccc} \rho_1(GAG^*) & \rho_2(GAG^*) &
\rho_3(GAG^*)
\\ \rho_2^*(GAG^*) & a_{22} & a_{23} \\ \rho_3^*(GAG^*) & a_{32} & a_{33} \end{array}
\right).
\end{eqnarray*}

Arguing as we did in the case that $rank(\psi)=2$, we see that
$\rho_2(D)=\rho_3(D)=0$ for all diagonal matrices $D \in M_2(\C)$,
so for $j=2,3$,
\begin{eqnarray}\label{rhoj}
\rho_j \left(\begin{array}{cc} b_{11} & b_{12} \\ b_{21} & b_{22}
\end{array} \right) = \rho_j\left(\begin{array}{cc} 0 & b_{12} \\ b_{21} &
0
\end{array} \right)\end{eqnarray} for all $B \in M_2(\C)$. For each $c$ on the unit circle $S^1$, let

\begin{equation*}w_{c} = \rho_2 \left(\begin{array}{cc} 0 & c \\ \bar{c}
& 0
\end{array} \right), \ \ z_{c}= \rho_3 \left(\begin{array}{cc} 0 & c \\ \bar{c}
& 0
\end{array} \right). \end{equation*}
Applying $\Phi$ to the family of positive $3 \times 3$ matrices
$\{M_c\}_{c \in S^1}$ defined by $M_c=e_{22} + c e_{23} + \bar{c}
e_{32}+ e_{33}$, we find
\begin{eqnarray*} 0 & \leq & det(\Phi(M_c)) =
-|w_c|^2-|z_c|^2 + 2 Re(c w_c \overline{z_c})
\\ & = & -|c w_c|^2-|z_c|^2+ 2 Re( c w_c \overline{z_c}) = - |c w_c-z_c|^2,
\end{eqnarray*}
hence $c w_c = z_c$ for all $c \in S^1$.  This gives us
\begin{eqnarray}\label{S1}  c^2 \rho_2\left(\begin{array}{cc}
0 & 1 \\ 0 & 0  \end{array} \right)  + |c|^2
\rho_2\left(\begin{array}{cc} 0 & 0 \\ 1 & 0  \end{array} \right)  =
c \rho_3\left(\begin{array}{cc} 0 & 1 \\ 0 & 0  \end{array} \right)
+ \bar{c} \rho_3\left(\begin{array}{cc} 0 & 0 \\ 1 & 0  \end{array}
\right)
\end{eqnarray}
for all $c \in S^1$.  Applying \eqref{S1} to $c= 1$ and $c=-1$ yields
\begin{eqnarray}\label{abc} \rho_j\left(\begin{array}{cc}
0 & 1 \\ 0 & 0  \end{array} \right)  = -
\rho_j\left(\begin{array}{cc} 0 & 0 \\ 1 & 0  \end{array} \right)
\end{eqnarray} for $j=2,3$.  Letting \begin{equation*} b=
\rho_2\left(\begin{array}{cc} 0 & 1 \\ 0 & 0  \end{array} \right) \
\textrm{ and } \ d=\rho_3\left(\begin{array}{cc} 0 & 1 \\ 0 & 0
\end{array} \right) ,\end{equation*} we rewrite \eqref{S1} as $$c^2
b - |c|^2 b = c d - \bar{c} d.$$ Applying this to $c=i$ and $c=-i$
yields $b=d=0$, whereby \begin{equation*}
\rho_j\left(\begin{array}{cc} 0 & 1 \\ 0 & 0  \end{array} \right) =0
\ \textrm{ and thus } \ \rho_j\left(\begin{array}{cc} 0 & 0 \\ 1 & 0
\end{array} \right) =0 \textrm{ by \eqref{abc} for } j=2,3.
\end{equation*}
We conclude from \eqref{rhoj} that $\rho_2=\rho_3 \equiv 0$, hence
$\Phi$ has the form
\begin{eqnarray*} \Phi(A)= \left(\begin{array}{ccc} \rho_1(GAG^*) & 0 &
0
\\ 0 & a_{22} & a_{23} \\ 0 & a_{32} & a_{33} \end{array}
\right).
\end{eqnarray*}
Since $\rho_1$ is a state on $M_2(\C)$, we know that for some
unitary $Y \in M_2(\C)$ and $\lambda \in [0,1]$, we have
\begin{equation*}
\rho_1 \Big(Y \left(
\begin{array}{cc} b_{11} & b_{12} \\ b_{21} & b_{22}
\end{array} \right) Y^* \Big) = \lambda b_{11} + (1 -\lambda)
b_{22} \end{equation*} for all $B \in M_2(\C)$, so for every $A \in
M_3(\C)$,

\begin{equation*}
\rho_1 \Big(Y [GAG^*] Y^* \Big) = \lambda a_{22} + (1 -\lambda)
a_{33}.
\end{equation*}
Letting
\begin{displaymath}
X = \left( \begin{array}{ccc}
1 & 0 & 0 \\ 0 &
\multicolumn{2}{c}{\multirow{2}{*}{$[Y]$}}
\\ 0
\end{array} \right),
\end{displaymath}
we observe that $\Phi_X$ has the form  \eqref{arg5}.

\qqed

\begin{lem}\label{annih}  Suppose
$\Phi: M_3(\C) \rightarrow M_3(\C)$ is a linear map which does
not annihilate any projections and satisfies $\rank(\Phi)>1$.
Then $\Phi \in \mathcal{E}_3$
if and only if, up to conjugacy, it has one of the
following forms for all $A \in M_3(\C)$:

\begin{equation} \label{arg10} \tag{IV}
\Phi(A) = \left( \begin{array}{ccc}
a_{11} & 0 & 0 \\ 0 & a_{22} & 0  \\ 0 & 0 & a_{33}
\end{array}
\right);
\end{equation}

\begin{equation} \label{arg7} \tag{V} \Phi(A) = \left( \begin{array}{ccc}
a_{11} & 0 & 0 \\ 0 & a_{22} & a_{23} \\ 0 & a_{32} & a_{33}
\end{array}
\right);
\end{equation}

\begin{equation} \label{arg8} \tag{VI} \Phi(A) = \left( \begin{array}{ccc} a_{11} & 0 & 0 \\ 0
& \lambda a_{22} + (1- \lambda) a_{33} & 0 \\
0 & 0 & \lambda a_{22} + (1- \lambda) a_{33} \end{array} \right), \ \ \lambda \in (0,1) ;
\end{equation}
\begin{equation} \label{arg9} \tag{VII} \Phi(A)=A. \end{equation}
\end{lem}

{\bf Proof:} The backward direction follows from inspection
of the maps \eqref{arg10} through \eqref{arg9}.  Assume
the hypotheses of the forward direction.
By Lemma \ref{fixordestroy},
$\Phi$ fixes a rank one projection $E$.  Note that $U^*EU = e_{11}
$ for some unitary $U \in M_3(\C)$, so
$$\Phi_U(e_{11}) = e_{11}.$$
Therefore, we may assume that $E=e_{11}$ and $\Phi(e_{11})=e_{11}$.
Let $F=I-E=e_{22}+e_{33}$.  For some functionals $\tau_2, \tau_3 \in
M_3(\C)^*$ and some linear map $\Psi: M_3(\C) \rightarrow M_2(\C)$,
we have

\begin{displaymath}
\Phi(A) = \left(\begin{array}{ccc} a_{11} & \tau_2(A) & \tau_3(A)
\\ \tau^*(A) & \multicolumn{2}{c}{\multirow{2}{*}{$\Big[\Psi(A)\Big]$}}
\\ \tau_3 ^*(A)
\end{array} \right).
\end{displaymath}

However, by Lemma \ref{form}, $\Phi$ satisfies
$$\Phi(A)  =  E\Phi(EAE)E +
E\Phi(EAF)F + F \Phi(FAE)E + F\Phi(FAF)F,$$ so
\begin{equation} \label{morerhoj}
\Psi(A)=\Psi(FAF), \ \ \tau_j(A)=\tau_j(EAF)=\tau_j(a_{12}e_{12} +
a_{13} e_{13})
\end{equation} for all $A \in M_3(\C)$ and $j=2,3$.  Let
\begin{equation*}
G=\left(\begin{array}{ccc} 0 & 1 & 0 \\
0 & 0 & 1  \end{array}\right) \in M_{2 \times 3}(\C) \ \textrm{ and } \ J=\left(\begin{array}{c} 1 \\
0
\\ 0     \end{array}\right) \in M_{3 \times 1}(\C).
\end{equation*}
Defining $\psi: M_2(\C) \rightarrow M_2(\C)$ and $\rho_j \in M_{1
\times 2}(\C)^*$ ($j=2,3$) by $\psi(B)=\Psi(G^*BG)$ and $\rho_j(C) =
\tau_j(JCG)$ for all $B \in M_2(\C)$ and $C \in M_{1 \times 2}(\C)$,
we see that $\Phi$ has the form
\begin{equation} \label{hah} \Phi(A) = \left(\begin{array}{ccc} a_{11} &
\rho_2\left(\begin{array}{cc} a_{12} & a_{13} \end{array} \right)
& \rho_3\left(\begin{array}{cc} a_{12} & a_{13} \end{array} \right) \\
\rho_2 ^* \left(\begin{array}{cc} a_{21} \\ a_{31} \end{array}
\right) &
\multicolumn{2}{c}{\multirow{2}{*}{$\Big[\psi\left(\begin{array}{cc}
a_{22} & a_{23} \\ a_{32} & a_{33} \end{array} \right)\Big]$}}
\\ \rho_3^*\left(\begin{array}{cc} a_{21} \\ a_{31} \end{array}
\right)
\end{array} \right)
\end{equation}
for all $A \in M_3(\C)$.

From equation \eqref{hah} and the fact that $\Phi^2=\Phi$, we have
$\psi^2=\psi$ and $\psi(I_2)=I_2$ for the $2 \times 2$ identity
matrix $I_2$. Moreover, $\psi$ is completely positive since
$\psi(B)= G \Phi(G^*BG)G^*$ for all $B \in M_2(\C)$.  Therefore,
$\psi \in \mathcal{E}_2$, so $\rank(\psi) \in \{1,2,4\}$ by
Proposition \ref{E2}.

{\bf Case (i):} If $\psi$ has rank one, then it has the form
\begin{displaymath}
\psi \left( \begin{array}{cc} a_{22} & a_{23}\\ a_{32} & a_{33}
\end{array} \right) = \rho
\left( \begin{array}{cc} a_{22} & a_{23}\\ a_{32} & a_{33}
\end{array} \right) I_2,
\end{displaymath} where $\rho$ is faithful since $\Phi$
does not annihilate any nonzero projections.  For all $A \in
M_3(\C)$ we have
\begin{displaymath}
\Phi(A)= \left(\begin{array}{ccc} a_{11} &
\rho_2\left(\begin{array}{cc} a_{12} & a_{13} \end{array} \right)
& \rho_3\left(\begin{array}{cc} a_{12} & a_{13} \end{array} \right) \\
\rho_2 ^* \left(\begin{array}{cc} a_{21} \\ a_{31} \end{array}
\right) & \multicolumn{2}{c}{\multirow{2}{*}{$\Big[\rho \left(
\begin{array}{cc} a_{22} & a_{23}\\ a_{32} & a_{33}
\end{array} \right) I_2\Big]$}}
\\ \rho_3^*\left(\begin{array}{cc} a_{21} \\ a_{31} \end{array}
\right)
\end{array} \right).
\end{displaymath}
Let $C$ be the matrix
\begin{displaymath}
C= \left( \begin{array}{cc} \rho_2(1 \ \ 0) & \rho_2(0 \ \ 1)
\\ \rho_3(1 \ \ 0) & \rho_3(0 \ \ 1) \end{array} \right).
\end{displaymath}
Since $\Phi^2(e_{12})=\Phi(e_{12})$ and $\Phi^2(e_{13}) =
\Phi(e_{13})$, we have $C^2=C$.

If $C=0$, then we repeat a familiar argument:  Since $\rho$ is
faithful and \begin{displaymath}\rho \left(  \begin{array}{cc}  1 &
0 \\ 0 & 1 \end{array} \right)=1 ,
\end{displaymath} we know that for some $2 \times 2$ unitary $T$ and $\lambda \in
(0,1)$,
\begin{equation*}
\rho \left( T \left( \begin{array}{cc}  a_{22} & a_{23}
\\ a_{32} & a_{33}
\end{array}\right) T^* \right)  = \lambda a_{22} + (1-\lambda) a_{33}
\end{equation*}
for all $A \in F M_3(\C) F$.  Letting
\begin{displaymath}
Z= \left( \begin{array}{ccc}
1 & 0 & 0
\\ 0 & \multicolumn{2}{c}{\multirow{2}{*}{$[T]$}}
\\ 0
\end{array} \right),
\end{displaymath}
we see that
\begin{equation*} \Phi_Z(A) = \left( \begin{array}{ccc} a_{11} & 0 & 0 \\ 0
& \lambda a_{22} + (1- \lambda) a_{33} & 0 \\
0 & 0 & \lambda a_{22} + (1- \lambda) a_{33}) \end{array} \right)
\end{equation*}
for all $A \in M_3(\C)$, so $\Phi_Z$ has the form \eqref{arg8}.

Now suppose $\rank(C) \geq 1$.  Since $C^2=C$, $C$ fixes a unit
vector $\vec{x}$,
\begin{equation*} \vec{x} = \left( \begin{array}{cc}
a \\ b  \end{array} \right), \ \ |a|^2+|b|^2=1.
\end{equation*}
%
%
%
In other words,
\begin{displaymath}
\left( \begin{array}{c} a \\ b \end{array} \right)
= C \left( \begin{array}{c} a \\ b \end{array} \right) =
\left( \begin{array}{c} \rho_2(a \ \ b) \\
\rho_3(a \ \ b) \end{array} \right).
\end{displaymath}
Letting
\begin{displaymath}A= \left(\begin{array}{c} 1 \\ \bar{a} \\ \bar{b}
\end{array} \right) \left(\begin{array}{ccc} 1 & a & b
\end{array} \right) =
\left( \begin{array}{ccc} 1 & a & b
\\ \bar{a} & |a|^2 & \bar{a}b
\\ \bar{b} & a \bar{b} & |b|^2
\end{array} \right) = \left( \begin{array}{ccc}
1 & a & b
\\ \bar{a} & \multicolumn{2}{c}{\multirow{2}{*}{$[P]$}}
\\ \bar{b} &
\end{array} \right) \geq 0,
\end{displaymath}
we have
\begin{displaymath}
\Phi(A)= \left( \begin{array}{ccc}
1 & a & b
\\ \bar{a} & \rho(P) & 0
\\ \bar{b} & 0 & \rho(P)
\end{array} \right) \geq 0,
\end{displaymath}
hence $0 \leq det(A) = \rho(P)(\rho(P) - |a|^2 - |b|^2)=
\rho(P)(\rho(P)-1)$. But $0 < \rho(P) \leq 1$ since $P$ is a rank
one projection, so $\rho(P)=1$. Therefore, $\Phi$ annihilates the
rank one projection
\begin{displaymath}
\left(\begin{array}{ccc} 0 & 0 & 0
\\ 0 &  \multicolumn{2}{c}{\multirow{2}{*}{$[I_{2} - P]$}}
\\ 0 \end{array} \right),
\end{displaymath}
contradicting our assumption that $\Phi$ does not destroy any
nonzero projections.

{\bf Case (ii):} If $\rank(\psi)=2$, then by Proposition \ref{E2},
$\psi_V$ is the diagonal map for some unitary $V \in M_2(\C)$.
Letting
\begin{displaymath}
S = \left( \begin{array}{ccc} 1 & 0 & 0 \\ 0 &
\multicolumn{2}{c}{\multirow{2}{*}{$[V]$}}
\\ 0
\end{array} \right),
\end{displaymath}
we see $\Phi_S(e_{11})=S^* \Phi(Se_{11}S^*)S=
S^*\Phi(e_{11})S^*=e_{11}$ and $\psi_V(B)=G \Phi_S(G^*BG)G^*$ for
all $B \in M_2(\C)$.  Since $\Phi_S$ fixes $e_{11}$ and does not
annihilate any nonzero projections, we may argue as we did earlier
in the proof (using Lemma \ref{form}) to conclude that for some
functionals $\rho_2'$ and $\rho_3'$ acting on $M_{1 \times 2}(\C)$,
$\Phi_S$ has the form
\begin{displaymath}
\Phi_S(A)=\left( \begin{array}{ccc} a_{11} &
\rho_2'\left(\begin{array}{cc} a_{12} & a_{13} \end{array}\right) &
\rho_3 '
\left(\begin{array}{cc} a_{12} & a_{13} \end{array}\right) \\
\rho_2'^*\left(\begin{array}{cc} a_{21} \\ a_{31} \end{array}\right)
& a_{22} & 0 \\ \rho_3'^*\left(\begin{array}{cc} a_{21} \\
a_{31} \end{array}\right) & 0 & a_{33}
\end{array} \right)
\end{displaymath}
Replacing $\Phi_S$ with $\Phi$ and erasing the primes from the
functionals $\rho_2'$ and $\rho_3'$, we continue our argument.

Now
\begin{displaymath}
\Phi \left(\begin{array}{ccc} 1 & 1 & 0
\\ 1 & 1 & 0 \\ 0 & 0 & 0
\end{array} \right) =
\left(\begin{array}{ccc} 1 & \rho_2(1 \ \ 0) & \rho_3(1 \ \ 0)
\\ \overline{\rho_2(1 \ \ 0)} & 1 & 0 \\ \overline{\rho_3(1 \ \ 0)} & 0 & 0
\end{array} \right)
\geq 0,
\end{displaymath}
hence $\rho_3(1 \ \ 0)=0$.  Similarly, positivity of $\Phi(e_{11}+
e_{13}+e_{31}+e_{33})$ implies $\rho_2(0 \ \ 1)=0$.  It follows that
for some $z_2, z_3 \in \C$, we have $\rho_2(a_{12} \ \ a_{13})= z_2
a_{12}$ and $\rho_3(a_{12} \ \ a_{13})= z_3 a_{13}$ for all $(a_{12}
\ \ a_{13}) \in M_{1 \times 2}(\C)$. Since $\Phi^2=\Phi$ we have
$z_j^2=z_j$, so $z_j \in \{0,1\}$ for $j=2,3$. Therefore, $\Phi$ is
the Schur map $\Phi(A)=M \bullet A$, where
\begin{displaymath}M=\left(\begin{array}{ccc} 1 & z_2 & z_3
\\ z_2 & 1 & 0 \\ z_3 & 0 & 1 \end{array}
\right) \geq 0.
\end{displaymath}
If $z_2=z_3=0$, then $\Phi$ has the form \eqref{arg10}.  If $z_2 =
1$, then by positivity of $M$ we have $z_3=0$, and we note that for
the unitary matrix $U= e_{13} + e_{22} + e_{31}$, $\Phi_U$ has the
form of \eqref{arg7}.  On the other hand, if $z_3=1$ then $z_2=0$ by
positivity of $M$.  Letting $V= e_{12} + e_{21} + e_{33}$, we
observe that $\Phi_V$ has the form of \eqref{arg7}.

{\bf Case (iii):} If $\psi$ is the identity map, we may repeat the
same argument we just used to show that for some $z_2, z_3 \in
\{0,1\}$, we have $\rho_2(a_{12} \ \ a_{13})= z_2 a_{12}$ and
$\rho_3(a_{12} \ \ a_{13})= z_3 a_{13}$ for all $(a_{12} \ \ a_{13})
\in M_{1 \times 2}(\C)$. Therefore, $\Phi(A) = N \bullet A$ for all
$A \in M_3(\C)$, where
\begin{displaymath}N=\left(\begin{array}{ccc} 1 & z_2 & z_3
\\ z_2 & 1 & 1 \\ z_3 & 1 & 1 \end{array}
\right) \geq 0.
\end{displaymath}
From positivity of $N$, we conclude that either $z_2=z_3=1$ (i.e.
$\Phi$ is the identity map \eqref{arg9}) or $z_2=z_3=0$ (in which
case $\Phi$ has the form \eqref{arg7}).

\qqed

Lemmas \ref{annih2} and \ref{annih} give us the following:
\begin{thm}\label{3limits}
A linear map $\Phi: M_3(\C) \rightarrow M_3(\C)$ is in
$\mathcal{E}_3$ if and only if $\Phi$ has the form
$\Phi(A)=\rho(A)I$ for some faithful state $\rho \in M_3(\C)^*$ or,
up to conjugacy, $\Phi$ is one of the maps \eqref{arg1} through
\eqref{arg9}.
\end{thm}

{\bf Proof:}  The only case not covered by Lemmas \ref{annih2} and
\ref{annih} is when $\Phi$ is a rank one map which does not
annihilate any nonzero projections.  It is clear that such a map
$\Phi$ is in $\mathcal{E}_3$ if and only if it is of the form
$\Phi(A)=\rho(A)I$ for a faithful state $\rho$. \qqed

\begin{cor}
Let $\phi: M_3(\C) \rightarrow M_3(\C)$ be unital and $q$-positive.
Then $\phi$ has rank $1, 2, 3, 4, 5,$ or $9$.
\end{cor}

\section{Classification of unital $q$-positive and $q$-pure maps on $M_2(\C)$}
If unital $q$-positive maps $\phi$ and $\psi$ are conjugate, then
their limits $L_\phi$ and $L_\psi$ are naturally conjugate as well:
\begin{lem}\label{Lphiu}
Suppose $\phi: M_n(\C) \rightarrow M_n(\C)$ is $q$-positive and
$||\phi(I + t \phi)^{-1}||<1$ for all $t>0$, and let $U \in M_n(\C)$ be unitary.
Then $L_{(\phi_U)} = (L_\phi)_U.$
\end{lem}
{\bf Proof:}  We know from Proposition 4.5 of \cite{Me} that
$$t \phi_U (I + t \phi_U)^{-1}(A) = U^* \phi(I + t \phi)^{-1}(UAU^*)U$$ for all $t>0$
and $A \in M_n(\C)$, so
\begin{eqnarray*}
L_{\phi_U}(A) & = & \lim_{t \rightarrow \infty} t \phi_U(I + t \phi_U)^{-1}(A)
= \lim_{t \rightarrow \infty} tU^* \phi(I + t \phi)^{-1}(UAU^*)U
\\ & = & U^*\Big[\lim_{t \rightarrow \infty} t\phi(I + t \phi)^{-1}(UAU^*)\Big] U =U^*L_\phi(UAU^*)U=(L_\phi)_U(A).
\end{eqnarray*}
\qqed

\begin{prop} \label{qpos} Let $\phi: M_2(\C) \rightarrow M_2(\C)$ be a unital
linear map of rank $2$.  Then $\phi$ is $q$-positive if and only if,
for some unitary $U \in M_2(\C)$ and numbers $\lambda \in (0,1], \lambda'
\in [0,1)$ with $\lambda> \lambda'$, we have
\begin{equation}\label{2x2eq} \phi_U(A) =
\left( \begin{array}{cc} \lambda a_{11} + (1-\lambda)a_{22} & 0 \\
0 & \lambda' a_{11} + (1-\lambda')a_{22} \end{array} \right)
\end{equation} for all $A=\sum_{i,j=1}^2 a_{ij}e_{ij} \in M_2(\C).$
\end{prop}

{\bf Proof:}  For the forward direction, assume that $\phi$ is $q$-positive.  It
follows from Proposition \ref{E2} that for some unitary $U \in
M_2(\C)$, $(L_{\phi})_U$ is the diagonal map
$$(L_{\phi})_U(A) = a_{11}e_{11} + a_{22}e_{22}.$$  But
$(L_\phi)_U = L_{(\phi_U)}$ by Lemma \ref{Lphiu} and
$\range(\phi_U) = \range(L_{\phi_U})$ by Lemma \ref{limit},
hence $\range(\phi_U) = \spann\{e_{11}, e_{22}\}$.  Therefore, for some
positive functionals $\rho_1, \rho_2 \in M_2(\C)^*$,
$$\phi_U(A)=\rho_1(A)e_{11} + \rho_2(A) e_{22}$$
for all $A \in M_2(\C)$, where $\rho_1$ and $\rho_2$ are states
since $\phi_U(I)=I$. Since
$\nullspace(\phi_U)=\nullspace(L_{\phi_U})$ by Lemma \ref{limit}, we
have $\phi_U(e_{12})=\phi_U(e_{21})=0$, so
$\rho_j(e_{12})=\rho_j(e_{21})=0$ for $j=1,2$.  Therefore, there are numbers
$\lambda, \lambda' \in [0,1]$ such that \begin{equation} \label{5}
\rho_1(A)= \lambda a_{11} + (1- \lambda)a_{22}, \ \ \rho_2(A)=
\lambda' a_{11} + (1- \lambda')a_{22}\end{equation} for all $A \in
M_2(\C)$. Let $Q=\lambda - \lambda'$, and for every $t \geq 0$, let
$$D_t= 1+t(1 + Q)+ t^2Q.$$  To prove the
forward direction, it suffices to show that $Q>0$, since it will
then automatically follow that $\lambda \in (0,1]$ and $\lambda' \in
[0,1)$.  For $j=1,2$, let $\nu_j \in M_2(\C)^*$ be
the functional $\nu_j(A)=a_{jj}$.  If $t \geq 0$ and $D_t \neq 0$,
then a straightforward computation shows that
$I + t \phi_U$ is invertible and
$$(I + t \phi_U)^{-1}(A)= A - \mu_{1,t}(A)e_{11} - \mu_{2,t}(A)e_{22}$$
for all $A \in M_2(\C)$, where $\mu_1$ and $\mu_2$ are the functionals
\begin{equation} \label{muse} \mu_{1,t} = \frac{t(\lambda +tQ)\nu_1+ t(1-\lambda)\nu_2}{D_t},  \  \  \mu_{2,t} =
\frac{t\lambda' \nu_1 + t(1-\lambda'+tQ)\nu_2}{D_t}.
\end{equation}
It follows that
\begin{equation}\label{ugh} t \phi_U (I + t \phi_U)^{-1}(A) = \mu_{1,t}(A)e_{11} + \mu_{2,t}(A)e_{22}
\end{equation}
for all $A \in M_2(\C)$.  If $Q=0$, then $\rank(\phi)=1$ by \eqref{5},
contradicting our assumption that $\rank(\phi)=2$. If $Q<0$, then
$D_{t_0}=0$ for some $t_0>0$. Since $||t \phi_U(I + t
\phi_U)^{-1}||<1$ for all $t>0$, the numerators of $\mu_{1,t}$ and
$\mu_{2,t}$ must both approach zero as $t \rightarrow t_0$.
With regard to $\mu_{1,t}$, this means that either $\lambda=1$ (contradicting
our assumption that $Q<0$) or
$$\nu_2 = - \frac{\lambda+t_0 Q}{1-\lambda} \nu_1,$$ which is clearly impossible.
Thus $Q>0$, proving the
forward direction.

Now assume the hypotheses of the backward direction.  For every $t>0$, we have $D_t
>0$, so $I + t \phi_U$ is invertible and $t \phi_U(I + t \phi_U)^{-1}$ has the form
\eqref{ugh}, where $\mu_{1,t}$ and $\mu_{2,t}$ are positive linear functionals
by \eqref{muse}.  Therefore, $\phi_U$ (and thus $\phi$) is
$q$-positive. \qqed

\begin{thm}\label{124}  Let $\phi: M_2(\C) \rightarrow M_2(\C)$ be a unital
linear map.  Then $\phi$ is $q$-positive if and only if it satisfies
one of the following:
\begin{enumerate}[(i)]
\item $\phi(A)=\rho(A)I$ for all $A \in M_2(\C)$, where $\rho \in M_2(\C)^*$
is a state.
\item For some $\lambda \in (0,1]$ and $\lambda' \in [0,1)$ with $\lambda> \lambda'$,
$\phi$ is conjugate to the map $\psi$ defined by
\begin{equation*}
\psi(A) =
\left( \begin{array}{cc} \lambda a_{11} + (1-\lambda)a_{22} & 0 \\
0 & \lambda' a_{11} + (1-\lambda')a_{22} \end{array} \right).
\end{equation*}
\item $\phi = \psi^{-1}$, where $\psi: M_2(\C) \rightarrow M_2(\C)$
is a unital conditionally negative map.
\end{enumerate}
\end{thm}

{\bf Proof:}  By Proposition \ref{E2}, we may assume that $\phi$ has
rank 1, 2, or 4.  From Proposition \ref{qpos} and Theorems \ref{statesbig}
and \ref{inverts},
conditions (i), (ii), and (iii) are the necessary and sufficient
conditions for $q$-positivity of unital linear maps $\phi: M_2(\C)
\rightarrow M_2(\C)$ of rank 1, 2, and 4, respectively. \qqed

Now that we have every unital $q$-positive $\phi: M_2(\C)
\rightarrow M_2(\C)$, we find that the only such maps which are
$q$-pure are rank one or invertible.
\begin{thm}\label{qpurem2}  A unital linear map $\phi: M_2(\C) \rightarrow M_2(\C)$ is
$q$-pure if and only if it satisfies one of the following:
\begin{enumerate}[(i)]
\item $\phi(A)=\rho(A)I$ for all $A \in M_2(\C)$, where $\rho \in
M_2(\C)^*$ is a faithful state;
\item  For some $\lambda \in \R$, $\phi$ is conjugate to the Schur map $\psi$ defined by
\begin{equation*} \psi(A) = \left(
\begin{array}{cc}
a_{11} & \frac{a_{12}}{1+ i \lambda}
\\ \frac{a_{21}}{1-i \lambda} & a_{22} \end{array}
\right)
\end{equation*}
for all $A \in M_2(\C)$.
\end{enumerate}
\end{thm}

{\bf Proof:}  By Theorems \ref{statesbig} and \ref{inverts},
conditions (i) and (ii) are the
necessary and sufficient conditions for a unital linear map of rank
1 or 4 to be $q$-pure.  Suppose that $\phi$ is a unital $q$-positive map of
rank 2, so by Theorem \ref{124}, it is conjugate to a map of the form \eqref{2x2eq}.
Since $q$-purity is invariant under conjugacy (Proposition 4.5 of \cite{Me}),
it suffices to assume
$\phi$ has the form \eqref{2x2eq} and show that $\phi$ is not $q$-pure.  Defining $\nu_1$
and $\nu_2$ as in the proof of Proposition \ref{qpos}, we recall that for every $t \geq 0$, we have
$t \phi(I + t \phi)^{-1}(A)= \mu_{1,t}(A)e_{11} +
\mu_{2,t}(A)e_{22},$ where $Q:=\lambda-\lambda'>0$ and

$$\mu_{1,t} = t \frac{(\lambda+tQ)\nu_1+ (1-\lambda)\nu_2}{1+ t(1+Q) +
t^2 Q} , \  \mu_{2,t} = t \frac{\lambda' \nu_1 + (1-\lambda'+tQ)\nu_2}
{1+ t(1+Q) + t^2 Q}.$$

Define $\Phi: M_2(\C) \rightarrow M_2(\C)$ by
$$\Phi(A) = \frac{Q \nu_1(A)}{1-\lambda'}e_{11}.$$  For every $t \geq 0$
and $A \in M_2(\C)$, we have $$(I+ t \Phi)^{-1}(A)= A - \frac{t Q \nu_1(A)}{1-\lambda'+tQ}
e_{11}$$ and
$$ \Phi(I + t \Phi)^{-1}(A) = \frac{ Q \nu_1(A)}{1-\lambda' + t Q}e_{11},$$
thus $\Phi \geq_q 0$.  Straightforward computations show
that $\phi-\Phi$ is completely positive and that for all $t>0$, we have
\begin{equation}\label{what} \Big(\phi(I + t\phi)^{-1}-\Phi(I + t \Phi)^{-1}\Big)(A)
= \eta_{1,t}(A)e_{11} + \frac{\mu_{2,t}(A)}{t}e_{22}\end{equation}
for all $A \in M_2(\C)$, where
$$\eta_{1,t}(A)  =  \frac{(\lambda+tQ)\nu_1(A) + (1-\lambda)\nu_2(A)}{1+t(1+Q)+t^2Q}
- \frac{Q \nu_1(A)}{1-\lambda'+tQ}.$$
Note that for every $t > 0$,
\begin{eqnarray*} \eta_{1,t} \geq 0 & \iff &
\frac{(\lambda+tQ)\nu_1 + (1-\lambda)\nu_2}{1+t(1+Q)+t^2Q}
\geq \frac{Q \nu_1}{1-\lambda'+tQ}
\\ & \iff & \frac{(\lambda +t(\lambda - \lambda'))\nu_1+ (1-\lambda) \nu_2}{1+
t(1+\lambda-\lambda') + t^2 (\lambda-\lambda')} \geq \frac{(\lambda-\lambda')\nu_1}{1-\lambda' + t(\lambda
-\lambda')}
\\ & \iff & (1-\lambda)\lambda' \nu_1
+ \Big((1-\lambda)(1-\lambda'+t(\lambda-\lambda')\Big)v_2 \geq 0.
\end{eqnarray*}
The coefficients of $\nu_1$ and $\nu_2$ are nonnegative in the above
expression, so $\eta_{1,t}$ is a positive linear
functional for all $t> 0$, hence $\phi \geq_q \Phi$ by \eqref{what}.
But $\rank(\Phi)=1$ while $\rank(\phi(I + s \phi)^{-1})=2$ for all $s \geq 0$,
so $\phi$ is not $q$-pure.
\qqed

\begin{prop}\label{notqpure} If $\phi: M_3(\C) \rightarrow M_3(\C)$ is a unital $q$-positive map
and $\phi(R)=0$ for some $R \gneq 0$, then $\phi$ is not $q$-pure.
\end{prop}
{\bf Proof:} If
$\phi(R)=0$ for some nonzero positive $R \in M_3(\C)$,  then $\phi$ annihilates a rank one projection $E$.
Letting $U \in
M_3(\C)$ be any unitary matrix such that $U^* E U = e_{11}$, we have $\phi_U(e_{11})=0$. Since $q$-purity is invariant
under conjugacy, we may replace $\phi_U$ with $\phi$ and continue
our argument.  Since $\phi(e_{11})=0$ we have $L_\phi(e_{11})=0$.
Replacing $L_\phi$ (and therefore $\phi$) with one of its conjugates
if necessary, we conclude $L_\phi$ has one of the forms \eqref{arg1}
through \eqref{arg5}.  Since $\phi$ and $L_\phi$ have the same range
and nullspace, it follows that \begin{equation} \label{duh}
\range(\phi) \subseteq \spann \{e_{11}, e_{22}, e_{23}, e_{32},
e_{33}\} \ \textrm{ and }  \ \phi(e_{1j})=\phi(e_{j1})=0 \textrm{
for } j=1,2,3.
 \end{equation}
Let $F=e_{22} + e_{33}$.  Line \eqref{duh} and Lemma \ref{FAF} imply
that for some state $\tau$ and some map $\Psi: M_3(\C) \rightarrow
M_2(\C)$, $\phi$ has the form

\begin{displaymath}
\phi (A) = \left( \begin{array}{ccc} \tau(A) & 0 & 0
\\ 0 & \multicolumn{2}{c}{\multirow{2}{*}{$[\Psi (A)]$}}
\\ 0 \end{array} \right),
\end{displaymath}
where $\tau(A)=\tau(FAF)$ and $\Psi(A)=\Psi(FAF)$ for all $A \in
M_3(\C)$.  Letting
\begin{displaymath}G = \left(\begin{array}{ccc}
0 & 1 & 0 \\ 0 & 0 & 1 \end{array}\right)
\end{displaymath}
and defining $\rho \in M_2(\C)^*$ and $\psi: M_2(\C) \rightarrow
M_2(\C)$ by $\rho(B)=\tau(G^*BG)$ and $\psi(B)=\Psi(G^*BG)$ for all
$B \in M_2(\C)$, we observe that $\phi$ has the form
\begin{displaymath}
\phi (A) = \left( \begin{array}{ccc} \rho(GAG^*) & 0 & 0
\\ 0 & \multicolumn{2}{c}{\multirow{2}{*}{$[\psi (GAG^*)]$}}
\\ 0 \end{array} \right).
\end{displaymath}
Note that $\psi$ has no negative eigenvalues. Indeed, suppose that
$\psi(B)= \lambda B$ for some $\lambda<0$ and $B \in M_2(\C)$. Let
$c= \rho(B)$. We see
$$\phi \Big( \frac{c}{\lambda} \  e_{11} + B \Big) = 0 + \phi(B)= ce_{11} + \lambda B
= \lambda \Big(\frac{c}{\lambda} \ e_{11} + B \Big),$$
contradicting the fact that $\phi$ has no negative eigenvalues.

Define $\phi': M_3(\C) \rightarrow M_3(\C)$ by $\phi'(A)=F\phi(A)F$
for all $A \in M_3(\C)$, so
\begin{displaymath}
\phi' (A) =
\left( \begin{array}{ccc} 0 & 0 & 0
\\ 0 & \multicolumn{2}{c}{\multirow{2}{*}{$[\psi (GAG^*)]$}}
\\ 0 \end{array} \right).
\end{displaymath}
We claim that $\phi'$ is $q$-positive.  Note that since $\psi$
has no negative eigenvalues, the same is true of $\phi'$.  Since
$\phi'$ commutes with $(I + t \phi')^{-1}$ for all $t \geq 0$, we have
\begin{eqnarray} \label{phi'} \phi'(I + t \phi')^{-1}(A) & = & (I + t \phi')^{-1}\phi'(A)=
(I + t \phi')^{-1}\phi'(FAF) \\ \nonumber & = & \phi'(I + t
\phi')^{-1}(FAF)\end{eqnarray} for all $A \in M_3(\C)$, and similarly,
$\phi(I + t \phi)^{-1}(A) = \phi(I + t \phi)^{-1} (FAF)$.

Let $A \in
M_3(\C)$. For some $2 \times 2$ matrix $B$, we have
\begin{equation} \label{hay}
(I + t \phi')^{-1}(FAF) = \left( \begin{array}{ccc} 0 & 0 & 0
\\ 0 &  \multicolumn{2}{c}{\multirow{2}{*}{$[B]$}}
\\ 0 &
\end{array} \right),  \ \
\phi'(I + t \phi')^{-1}(FAF) = \left( \begin{array}{ccc} 0 & 0 & 0
\\ 0 &  \multicolumn{2}{c}{\multirow{2}{*}{$[\psi(B)]$}}
\\ 0 &
\end{array} \right)
\end{equation}
while for some $d \in \C$, we have
\begin{equation} \label{hay2}
(I + t \phi)^{-1}(FAF) = \left( \begin{array}{ccc} d & 0 & 0
\\ 0 &  \multicolumn{2}{c}{\multirow{2}{*}{$[B]$}}
\\ 0 &
\end{array} \right), \ \
\phi(I + t \phi)^{-1}(FAF) = \left( \begin{array}{ccc} \rho(B) & 0 & 0
\\ 0 &  \multicolumn{2}{c}{\multirow{2}{*}{$[\psi(B)]$}}
\\ 0 &
\end{array} \right).
\end{equation}
Combining \eqref{phi'}, \eqref{hay}, and \eqref{hay2}, we find that
for all $A \in M_3(\C)$ and $t \geq 0$:
\begin{eqnarray} \label{toomany} \phi'(I + t \phi')^{-1}(A)
& = & \phi'(I + t \phi')^{-1}(FAF)
= F \Big( \phi(I + t \phi)^{-1}(FAF) \Big) F \\ \nonumber
&=& F \Big( \phi(I + t \phi)^{-1}(A) \Big) F. \end{eqnarray}
This shows that $\phi'$ is $q$-positive.  Furthermore, from equation \eqref{toomany}
and the fact that $$\phi(I + t \phi)^{-1}(A)
= e_{11} \Big( \phi(I + t \phi)^{-1}(A)\Big) e_{11} +
F \Big( \phi(I + t \phi)^{-1}(A) \Big)F$$
for all $A \in M_3(\C)$, we find that
$$
\phi(I + t \phi)^{-1}(A) - \phi'(I + t \phi')^{-1}(A)  =
e_{11} \Big( \phi(I + t \phi)^{-1}(FAF) \Big) e_{11}.
$$
Since the last line is the composition of completely positive maps
for every $t \geq 0$, we have $\phi \geq_q \phi'$.  Finally, we note
that $e_{11}  \phi'(I)  e_{11} = 0$, whereas for every $s \geq 0$,
$$e_{11} \Big( \phi(I + s \phi)^{-1}(I)\Big) e_{11} =
e_{11} \Big(\frac{1}{1+s} \ I\Big) e_{11} = \frac{1}{1+s} \
e_{11}.$$ Therefore, $\phi'$ is not equal to $\phi(I + s \phi)^{-1}$
for any $s \geq 0$, so $\phi$ is not $q$-pure. \qqed

\section{A cocycle conjugacy result}
Let $\nu$ be a type II Powers weight of the
form $\nu(\sqrt{I - \Lambda(1)}B \sqrt{I - \Lambda(1)}) = (f,Bf)$.
Suppose $\phi: M_n(\C) \rightarrow M_n(\C)$ ($n \geq 2$) and $\psi:
M_k(\C) \rightarrow M_k(\C)$ are unital and $q$-positive, where
$\rank(\phi)=1$ and $\psi$ is invertible.  We have seen that if
$\phi$ and $\psi$ are $q$-pure, then they are fundamentally
``different'' in the sense that $(\phi, \nu)$ and $(\psi, \nu)$ induce
non-cocycle conjugate $E_0$-semigroups (a consequence of Theorems \ref{statesbig} and \ref{inverts}).  We
now find that the previous sentence holds if we remove the assumption that $\phi$ and
$\psi$ are $q$-pure.  In fact, we may replace the assumption that
$\psi$ is invertible with the much weaker assumption that $L_\psi$
is a Schur map:

\begin{thm}\label{genl}
Let $\phi: M_n(\C) \rightarrow M_n(\C)$ ($n \geq 2$) and $\psi: M_k(\C)
\rightarrow M_k(\C)$ be unital $q$-positive maps.  Suppose that $\phi$
has rank one and that $L_\psi$ is a Schur map.  Let $\nu$ be a type II Powers weight
of the form $\nu(\sqrt{I - \Lambda(1)}
B \sqrt{I - \Lambda(1)}) = (f,Bf)$.

Then $(\phi, \nu)$ and $(\psi, \nu)$ induce non-cocycle conjugate $E_0$-semigroups.
\end{thm}

{\bf Proof:}  Let
$\alpha^d$ and $\beta^d$ be the $E_0$-semigroups induced by $(\phi, \nu)$
and $(\psi, \nu)$, respectively.  Suppose there is a nonzero $q$-corner $\gamma$ from $\phi$ to $\psi$, so $\Theta$
below is $q$-positive:
\begin{displaymath} \Theta= \left( \begin{array}{cc} \phi & \gamma
\\ \gamma^* & \psi \end{array} \right).
\end{displaymath}
Note that
\begin{displaymath} L_\Theta= \left( \begin{array}{cc} \phi & \sigma
\\ \sigma^* & L_\psi \end{array} \right),
\end{displaymath}
where $\sigma = \lim_{t \rightarrow \infty} t \gamma(I + t
\gamma)^{-1}$ is a corner from $\phi$ to $L_\psi$ since $L_\Theta$
is completely positive.  Furthermore, $\sigma^2= \sigma$, $\gamma=
\sigma \circ \gamma = \gamma \circ \sigma$, $\range(\sigma) =
\range(\gamma)$, and $\nullspace(\sigma)=\nullspace(\gamma)$.

Of course, $\phi$ has the form $\phi(A)=\rho(A)I$ for some state
$\rho \in M_n(\C)^*$.  Suppose that $\rho$ is faithful.  Let
$A \in M_{n \times k}(\C)$
be any norm one matrix in the range of $\sigma$, and
let $P \in M_n(\C)$ be the orthogonal
projection onto $\range(A) \subseteq \C^n$, so
$PA=A$ and $A^*P=A^*$.  Applying
$L_\Theta$ to the positive matrix $Q \in M_{n+k}(\C)$ given by
\begin{displaymath} Q = \left( \begin{array}{cc} P & 0
\\ 0 & I_{k}  \end{array} \right)
\left( \begin{array}{cc} I_n & A
\\ A^* & I_{k} \end{array} \right)
\left( \begin{array}{cc} P & 0
\\ 0 & I_{k}  \end{array} \right) =
\left( \begin{array}{cc} P & PA
\\ A^*P & I_{k}  \end{array} \right)=
\left( \begin{array}{cc} P & A
\\ A^* & I_{k}  \end{array} \right),
\end{displaymath}
we see from complete positivity of $L_\Theta$ that
\begin{displaymath} L_\Theta(Q)=
\left( \begin{array}{cc} \rho(P)I & A
\\ A^* & I_{k} \end{array} \right) \geq 0.
\end{displaymath}
Since $||A||=1$, positivity of the above matrix implies that
$\rho(P)=1$, hence $P=I_n$ by faithfulness of $\rho$.  Since $P$
is the orthogonal projection onto the range of $A$, we have
$\rank(A)=n$.  We conclude that every nonzero element
of $\range(\sigma)$ has rank $n$.

For some matrix unit $e_{ij} \in M_{n \times k}(\C)$,
we have $M:= \sigma(e_{ij}) \neq 0$, so $\rank(M)=n$.
By complete positivity of $L_\Theta$, the matrix $R$
below must be positive:
\begin{displaymath}
R = L_\Theta \left( \begin{array}{cc} e_{ii} & e_{ij}
\\ e_{ji} & e_{jj} \end{array} \right)
= \left( \begin{array}{cc} \rho(e_{ii})I_n & M
\\ M^* & e_{jj} \end{array} \right).
\end{displaymath}
However, $R$ is not positive.  Indeed, since
$\rank(M) = n \geq 2$, there exists a vector $g \in \C^k$ such
that $e_{jj}g=0$ but $M g \neq 0$.  For all $\lambda \in \R$
, we have
\begin{displaymath}\Big\langle \left( \begin{array}{c} Mg \\ -\lambda g \end{array} \right),
 R
\left( \begin{array}{c} Mg \\ -\lambda g \end{array} \right) \Big\rangle
=\Big( \rho(e_{ii}) - 2 \lambda\Big) ||Mg||^2,
\end{displaymath}
which is negative whenever $\lambda > 1$.  We conclude
$R \ngeq 0$, contradicting complete positivity of $L_\Theta$.  Therefore,
there is no nonzero $q$-corner from $\phi$ to $\psi$, so $\alpha^d$
and $\beta^d$ are non-cocycle conjugate by Proposition \ref{hypqc}.

Now suppose that $\rho$ is not faithful, so for some mutually
orthogonal norm one vectors $\{f_i\}_{i=1}^{p<n} \subset \C^n$
and positive numbers $\lambda_1, \ldots, \lambda_p$ with $\sum_{i=1}^p \lambda_i =1$, we have
$\rho(A) = \sum_{i=1}^{p} \lambda_i (f_i, A f_i)$ for all $A \in M_n(\C)$.
For some unitary $U \in M_n(\C)$ we have $$\phi_U(A)=\rho(UAU^*)I =
\Big(\sum_{i= n-p+1}^{n} \lambda_{i-n+p} \  a_{ii} \Big) I.$$
Since the $E_0$-semigroup $\alpha_U ^d$ induced by
$(\phi_U, \nu)$ is cocycle conjugate to $\alpha^d$ by Proposition \ref{dumb}, the theorem follows
if we show that $\alpha_U^d$ is not cocycle conjugate to $\beta^d$.

If there is a hyper maximal $q$-corner $\gamma$ from $\phi_U$ to $\psi$,
then
\begin{displaymath} \Theta= \left( \begin{array}{cc} \phi_U & \gamma
\\ \gamma^* & \psi \end{array} \right) \geq_q 0,
\end{displaymath}
and we have
\begin{displaymath} L_\Theta= \left( \begin{array}{cc} \phi_U & \sigma
\\ \sigma^* & L_\psi \end{array} \right),
\end{displaymath}
where $\sigma = \lim_{t \rightarrow \infty} t \gamma(I + t
\gamma)^{-1}$ is a norm one corner from $\phi_U$ to $L_\psi$ such
that $\sigma^2= \sigma$ and $\range(\sigma) = \range(\gamma)$.

Note that $\phi_U(e_{11}) = 0$, hence
$L_\Theta(e_{11})=0$, so by Lemma \ref{FAF} we have
\begin{displaymath}
\sigma \left( \begin{array}{cccc} b_{11} & b_{12} & \cdots & b_{1k}
\\ 0 & 0 & \cdots & 0
\\ \vdots & \vdots & \vdots & \vdots
\\ 0 & 0 & \cdots & 0 \end{array} \right) \equiv 0.
\end{displaymath}

Therefore, for some $\ell: M_{(n-1) \times k}(\C) \rightarrow M_{1 \times k}(\C)$
and $\tilde{\sigma}: M_{(n-1) \times k}(\C) \rightarrow M_{(n-1) \times k}(\C)$,
we have
\begin{displaymath}
\sigma \left( \begin{array}{c} B_{1 \times k}
\\ A_{(n-1) \times k} \end{array} \right)
= \left( \begin{array}{c} \ell(A_{(n-1) \times k})
\\ \tilde{\sigma}(A_{(n-1) \times k}) \end{array} \right).
\end{displaymath}
Note that since $\sigma^2=\sigma$ and $||\sigma||=1$,
we have $\tilde{\sigma}^2=\tilde{\sigma}$ and $||\tilde{\sigma}||=1$.

We claim that $\ell \equiv 0$.  To show this, we let $2 \leq i \leq n$
and $1 \leq j \leq k$ be arbitrary.  Since $L_\Theta$ is completely positive,
we have
\begin{displaymath}
0 \leq R:= L_\Theta \left( \begin{array}{cc} e_{ii} & e_{ij}
\\ e_{ji} & e_{jj} \end{array} \right)
= \left( \begin{array}{cc} \rho(e_{ii}) I_n & \sigma(e_{ij})
\\ \sigma^*(e_{ji}) & e_{jj} \end{array} \right).
\end{displaymath}
A rank argument similar to the one from the faithful case shows that
$\rank(\sigma(e_{ij}))\leq 1$ since $R$ is positive. If
$\sigma(e_{ij})=0$, then $\ell(e_{ij})=0$.  If
$\rank(\sigma(e_{ij}))=1$, we see from the form of $R$ that
$\sigma(e_{ij})$ is a column matrix of the form below for some
scalars $c_1, \ldots, c_n$:
$$\sigma(e_{ij}) = \sum_{m=1}^n c_m e_{mj}.$$
Since $\sigma^2=\sigma$ and $\sigma(e_{1j})=0$, we have
\begin{equation} \label{inc} \sum_{m=1}^n c_m e_{mj} = \sigma\Big(\sum_{m=1}^n c_m
e_{mj}\Big) = \sigma\Big(\sum_{m=2}^n c_m e_{mj}\Big).\end{equation}

If $c_1 \neq 0$, then by equation \eqref{inc},
$$ \Big|\Big|\sigma\Big(\sum_{m=2}^n c_m e_{mj}\Big)\Big|\Big| =\Big|\Big|\sum_{m=1}^n c_m e_{mj}
\Big|\Big|
>\Big|\Big|\sum_{m=2}^n c_m e_{mj}\Big|\Big|,$$
contradicting the fact that $\sigma$ is a contraction.  Hence
$c_1=0$, that is, $\ell(e_{ij})=0$.  Since $i$ and $j$ were chosen
arbitrarily, we conclude $\ell \equiv 0$.  This means that
\begin{equation} \label{alm}
\range(\sigma) \ \bigcap \ \spann \{e_{11}, e_{12}, \ldots, e_{1n} \} = \{0\}.
\end{equation}
The same holds for $\gamma$ since $\range(\gamma)=\range(\sigma)$.
Define $\Theta': M_{n + k}(\C) \rightarrow M_{n+k}(\C)$ by
$$\Theta'(A)= (I_{n+k} - e_{11}) \Theta(A) (I_{n+k} - e_{11}).$$
From equation \eqref{alm}, 
we have
\begin{displaymath} \Theta' = \left( \begin{array}{cc} \phi' & \gamma
\\ \gamma^* & \psi \end{array} \right),
\end{displaymath}
where $\phi': M_n(\C) \rightarrow M_n(\C)$ is the map
$\phi'(A)=\rho(A)(I_n-e_{11})$. Note $(\phi')^2= \phi'$ and $\phi'$
is $q$-positive, so $I + t \Theta'$ is invertible for all $t \geq 0$
and
$$\Theta'(I + t \Theta')^{-1} (A) = (I_{n+k} - e_{11})
\Big[ \Theta(I + t \Theta)^{-1}(A) \Big] (I_{n+k} - e_{11})$$
for all $A \in M_{n+k}(\C)$.  This shows that $\Theta' \geq_q 0$.  Furthermore,
$\Theta \geq_q \Theta'$ since $$\Theta'(I + t \Theta')^{-1} -
\Theta(I + t \Theta)^{-1}(A) = e_{11} \Theta(I + t \Theta)^{-1}(A) e_{11}$$
for all $t \geq 0$, $A \in M_{n+k}(\C)$.  Trivially, $\phi' \neq \phi$,
contradicting hyper maximality of $\gamma$.  We conclude there is no hyper maximal
$q$-corner from $\phi_U$ to $\psi$, hence $\alpha^d$ and $\beta^d$ are non-cocycle conjugate
by Proposition \ref{hypqc}.

\qqed
We conclude with the following:
\begin{cor} \label{lessgen}
Let $\phi_1: M_2(\C) \rightarrow M_2(\C)$ be a unital rank one
$q$-positive map.   Let $\phi_2: M_2(\C) \rightarrow M_2(\C)$ be the diagonal map,
and let $\phi_3: M_2(\C) \rightarrow M_2(\C)$ be a unital invertible $q$-positive
Schur map.  Suppose $\nu$ is a type II Powers weight of the form
$\nu(\sqrt{I - \Lambda(1)}
B \sqrt{I - \Lambda(1)}) = (f,Bf)$.

The boundary weight doubles $(\phi_i, \nu)$ and $(\phi_j, \nu)$ induce cocycle conjugate
$E_0$-semigroups if and only if $i=j$.
\end{cor}
{\bf Proof:}  For each $i$, let $\alpha_i ^d$ be the $E_0$-semigroup induced by
$(\phi_i, \nu)$.  Theorem \ref{genl} implies that $\alpha_1^d$ is not cocycle
conjugate to $\alpha_2^d$ or $\alpha_3^d$.  A result at the end
of \cite{bigpaper} 
shows that 
$\alpha_2 ^d$ and $\alpha_3 ^d$ are non-cocycle conjugate, but we present a proof
here for the sake of completeness.
Let $\gamma$ be any $q$-corner from $\phi_2$ to $\phi_3$, so $\Theta: M_4(\C) \rightarrow
M_4(\C)$ below is $q$-positive:
\begin{displaymath}
\Theta = \left( \begin{array}{cc} \phi_2 & \gamma \\ \gamma^* &
\phi_3   \end{array} \right).
\end{displaymath}
Applying $\Theta$ to the matrices $e_{11} + e_{1j} + e_{j1} + e_{jj}$
and $e_{22} + e_{2k} + e_{k2} + e_{kk}$ for $j,k=3,4$, we conclude from
completely positivity of
$\Theta$ that $\gamma$ is a Schur map.

Form $L_\Theta$, observing that
\begin{displaymath}
L_\Theta = \left( \begin{array}{cc} \phi_2 & \sigma \\ \sigma^* &
Id_{2 \times 2}   \end{array} \right),
\end{displaymath}
where $\sigma=\lim_{t \rightarrow \infty} t \gamma(I + t
\gamma)^{-1}$, $\sigma^2=\sigma$, $\range(\sigma)=\range(\gamma)$,
and $\nullspace(\sigma)=\nullspace(\gamma)$. Note that $\sigma$ is
also a Schur map, and write $\sigma$ in the form
\begin{displaymath}
\sigma(A) = \left( \begin{array}{cc} z_{11} a_{11} & z_{12}a_{12} \\ z_{21}a_{21} &
z_{22}a_{22}   \end{array} \right).
\end{displaymath}  Since $\sigma^2=\sigma$, we have $z_{ij} \in \{0,1\}$ for each $i$ and $j$.

We claim that
\begin{equation} \label{zzz} z_{21}=z_{22}=1 \ \textrm{ or } \ \ z_{21}=z_{22}=0. \end{equation}
To prove this, first suppose that
$z_{21}=1$.  Let $T \in M_4(\C)$
be the positive matrix whose entries are all $1$. Then
\begin{displaymath}
0 \leq L_\Theta \Big((I-e_{11})T(I-e_{11})\Big) = \left( \begin{array}{cccc}
 0 & 0 & 0 & 0 \\ 0 & 1 & 1 & z_{22} \\ 0 & 1 & 1 & 1
\\ 0 & z_{22} & 1 & 1 \end{array} \right).
\end{displaymath}
Since the above matrix is positive, the determinant of its bottom right $3 \times 3$ minor
must be nonnegative.  This quantity is $-(z_{22}-1)^2$, hence $z_{22}=1$.  On the other
hand, if $z_{21}=0$, then
$$0 \leq \det\Big[L_\Theta\Big((I - e_{11})T(I - e_{11})\Big)\Big] = -(z_{22})^2,$$ so
$z_{22}=0$, yielding \eqref{zzz}.

Analogous observations regarding $L_\Theta\Big((I - e_{22})T(I - e_{22})\Big)$  show that
\begin{equation}\label{zzz2} z_{11}=z_{12}=1 \ \textrm{ or } \ z_{11}=z_{12}=0.
\end{equation}
By equations \eqref{zzz} and \eqref{zzz2}, $\sigma$ is the Schur mapping
$\sigma(A)=M_j \bullet A$ for one of the three matrices below:
\begin{displaymath}
M_1 = \left( \begin{array}{cc} 1 & 1 \\ 1 & 1  \end{array} \right), \ \
M_2 = \left( \begin{array}{cc} 1 & 1 \\ 0 & 0  \end{array} \right), \
\textrm{ or }
M_3 = \left( \begin{array}{cc} 0 & 0 \\ 1 & 1  \end{array} \right).
\end{displaymath}

Note that $\sigma(A)=M_1 \bullet A$ is not a corner from $\phi_2$ to $Id_{2 \times 2}$,
since in that case we would have
\begin{displaymath} L_\Theta(T) = \left(
\begin{array}{cccc} 1 & 0 & 1 & 1 \\ 0 & 1 & 1 & 1 \\ 1 & 1 & 1 & 1
\\ 1 & 1 & 1 & 1 \end{array} \right) \ngeq 0,
\end{displaymath}
contradicting complete positivity of $L_\Theta$.

Recall that $\gamma$ is a Schur map and $\range(\sigma)
= \range(\gamma)$.  Therefore, if $\sigma(A) = M_2 \bullet A$ for all $A \in M_2(\C)$, then
$\gamma(A)=R \bullet A$ for some matrix $R$ such that $r_{21}=r_{22}=0$.
Letting $S=e_{11} + e_{33} + e_{44}$ and defining $\Theta'$ by
$$\Theta'(A)=S \Theta(A) S$$ for all $A \in M_4(\C)$, we see that $\Theta'$ is completely positive by
construction and
\begin{displaymath} \Theta'= \left( \begin{array}{cc} \phi_2' & \gamma
\\ \gamma^* & \phi_3, \end{array} \right)
\end{displaymath}
where $\phi_2'(A) = a_{11} e_{11}$ for all $A \in M_2(\C)$.  Furthermore, we have
$$\Theta'(I + t \Theta')^{-1}(A) = S \Big( \Theta(I + t \Theta)^{-1}(A) \Big)S$$
for all $t >0$ and $A \in M_4(\C)$, so $\Theta'$ is $q$-positive.  Also, $\Theta \geq_q \Theta'$
since $$\Theta(I + t \Theta)^{-1}(A) - \Theta'(I + t \Theta')^{-1}(A)= \frac{1}{1+t} e_{22}Ae_{22}$$ for all
$A \in M_4(\C)$ and $t \geq 0$.  Therefore, $\gamma$ is not a hyper maximal $q$-corner.

If $\sigma(A) = M_3 \bullet A$, then we argue precisely as we just did, noting first that
$\gamma$ is a Schur map $\gamma(A)= Y \bullet A$ for some $Y \in M_2(\C)$ with $y_{11}=
y_{12}=0$.   Letting $X=e_{22} + e_{33} + e_{44}$ and defining $\phi_2''$ by $\phi_2''(A)
= a_{22}e_{22}$, we note that the map $\Theta'': M_4(\C) \rightarrow M_4(\C)$ defined by
\begin{displaymath}
\Theta''= \left( \begin{array}{cc} \phi_2'' & \gamma \\ \gamma^* & \phi_3 \end{array} \right)
\end{displaymath}
satisfies $\Theta \neq \Theta''$ and $\Theta \geq_q \Theta''$.  We conclude that
$\gamma$ is not a hyper maximal $q$-corner.

We have shown that no $q$-corner $\gamma$ from $\phi_2$ and $\phi_3$ is hyper maximal, hence
$\alpha_2^d$ and $\alpha_3^d$ are non-cocycle conjugate by Proposition \ref{hypqc}.   \qqed

\begin{center}{\bf Acknowledgments}\end{center}
The author would like to thank Robert Powers for his guidance and
enthusiastic interest in this research.


\begin{thebibliography}{99}
\bibitem{arveson} W.B. Arveson, \emph{The Index of a Quantum Dynamical
Semigroup}, J. Funct. Anal. {\bf 146} (1997), 557-588.
\bibitem{arvindex} W.B. Arveson, \emph{Continuous Analogues
of Fock space}, Memoirs Amer. Math. Soc. {\bf 80}, no. 409 (1989).
\bibitem{Bhat}  B.V.R. Bhat, \emph{An index theory for quantum dynamical
semigroups}, Trans. A.M.S. {\bf 348} (1996), no. 2, 561-583.
\bibitem{choi}  M. Choi,
\emph{Completely positive linear maps on complex matrices}, Lin.
Alg. Appl. {\bf 10} (1975), 285-290.
\bibitem{izumi} M. Izumi, \emph{A perturbation problem for the shift
semigroup}, J. Funct. Anal. {\bf 251} (2007), 498-545.
\bibitem{DJ} C. Jankowski and D. Markiewicz, \emph{Gauge groups of certain
$E_0$-semigroups obtained from boundary weight doubles}, in
preparation.
\bibitem{Me} C. Jankowski, \emph{On type II$_0$ $E_0$-semigroups
induced by boundary weight doubles}, J. Func. Anal. {\bf 258}
(2010), no. 10, 3413-3451.
\bibitem{markie} D. Markiewicz and R.T. Powers, \emph{Local
unitary cocycles of $E_0$-semigroups}, J. Funct. Anal. {\bf 256}
(2009), no. 5, 1511-1543.
\bibitem{paulsen}  V.I. Paulsen,  \emph{Completely bounded maps and
dilations}, Longman Scientific and Technical, Essex, England, 1986.
\bibitem{hugepaper} R.T. Powers, \emph{Continous spatial semigroups of
completely positive maps of $B(H)$}, New York J. Math.
{\bf 9} (2003), 165-269.
\bibitem{bigpaper}  R.T. Powers, \emph{Construction of $E_0$-semigroups of
$B(\mathfrak{H})$ from $CP$-flows}, Advances
in Quantum Dynamics, Contemp. Math. {\bf 335}, Amer. Math. Soc., Providence, RI
(2003), 57-97.
\bibitem{T2}  B. Tsirelson, \emph{Non-isomorphic product systems}, Advances
in Quantum Dynamics, Contemp. Math. {\bf 335}, Amer. Math. Soc., Providence, RI
(2003), 273-328.
\bibitem{wigner}  E.P. Wigner, \emph{On unitary representations
of the inhomogeneous Lorentz group,} Ann. of Math. {\bf 40} (1939),
149-204.
\end{thebibliography}
\end{document}